\documentclass[%
 aip,
sd,%
 amsmath,amssymb,
reprint,%
]{revtex4-1}
\pdfoutput=1
\usepackage{graphicx}
\usepackage{dcolumn}
\usepackage{bm}
\usepackage{epstopdf}
\usepackage{subcaption}
\usepackage{color}
\newcommand{\change}[1]{{\color{black}#1}}
\newcommand{\R}{\mathbb{R}}
\DeclareMathOperator{\sech}{sech}

\draft 

\begin{document}


\title{Early-warning indicators for rate-induced tipping} 



\author{Paul Ritchie}
\email[]{pdlr201@exeter.ac.uk}

\author{Jan Sieber}
\email[]{J.Sieber@exeter.ac.uk}
\affiliation{Centre for Systems, Dynamics and Control, College of Engineering, Mathematics and Physical Sciences, Harrison Building, University of Exeter, Exeter, EX4 4QF, United Kingdom}


\date{\today}

\begin{abstract}
  A dynamical system is said to undergo rate-induced tipping when
    it fails to track its quasi-equilibrium state due to an
    above-critical-rate change of system parameters.  We study a
  prototypical model for rate-induced tipping, the saddle-node normal
  form subject to time-varying equilibrium drift and
  noise. We find that both most commonly used early-warning
    indicators, increase in variance and increase in autocorrelation,
    occur not when the equilibrium drift is fastest but with a
    delay. We explain this delay by demonstrating that the most likely
    trajectory for tipping also crosses the tipping threshold with a
    delay and therefore the tipping itself is delayed. We find solutions of the variational problem determining
  the most likely tipping path using numerical continuation
  techniques.  The result is a systematic study of the \change{most likely tipping time}
  in the plane of two parameters, distance from tipping threshold and noise
  intensity.
\end{abstract}


\maketitle 

\begin{quotation}
  The notion of tipping describes the phenomenon that at certain
  critical levels or rates the output of a system changes
  disproportionately compared to the change in input. Two particular
  examples in climate science are the possible collapse of the \change{Atlantic Meridional Overturning Circulation (AMOC)} due to increasing freshwater input or the sudden release of
  carbon in peatlands due to an external temperature increase above a
  critical rate (the compost bomb instability
  \cite{wieczorek2011excitability}).  There is an ongoing debate, for
  example, in climate science \cite{boulton2014early} and
  ecology \cite{xu2015local} whether it is possible to find early-warning
  indicators robustly in time series of system outputs
  \emph{before} tipping occurs. 

  This paper focuses on the case of rate-induced tipping, which
  describes the scenario where a system fails to track its equilibrium
  due to a rapid change in parameters.  We show that two popular
  candidates for early-warning indicators in time series, an increase
  in autocorrelation and an increase in variance, appear to give a
  delayed warning signal for tipping. We study the phenomenon by
  looking at the interaction of rate-induced tipping and noise. We
  find that the most likely time for the noise to kick the system over
  the threshold (which is a curve in phase space) is \emph{after} the
  point where the threshold is closest to the equilibrium. We
  investigate this tipping (and, correspondingly, early-warning) delay
  systematically depending on two parameters: the distance of the
  parameter drift speed from its critical value and the noise
  intensity. We find that the delay is larger for smaller
  noise-intensity.
\end{quotation}

\section{Introduction}

Tipping events are often described as sudden, disproportionate changes
in output levels caused by small changes to input levels
\cite{ashwin2015parameter}. These can be irreversible events that have
huge, unwanted consequences. Therefore, the study of early-warning
  indicators is of great interest and so recent research has
developed and analyzed early-warning indicators, see
\citet{lenton2011early} for a review up to 2011 and
  \citet{williamson2015detection} for references to later results. A few currently
  debated examples of complex systems deemed vulnerable to tipping
  from climate science, ecology and financial markets are: the abrupt
reductions in Arctic summer sea ice \cite{holland2006future}, the
collapse of the Atlantic Meridional Overturning Circulation
(AMOC)
\cite{held2004detection}, the dieback of the Amazon rainforest
\cite{malhi2009exploring}, Australian ecosystems\cite{laurance201110},
the light-driven regime shifts in polar ecosystems
\cite{clark2013light}, the collapse of coral reefs due to global
warming and ocean acidification \cite{hoegh2007coral}, and crashes and
rebounds of financial markets \cite{yan2010diagnosis}. See also
\citet{lenton2008tipping} for a list of policy-relevant tipping
elements in the climate system. \citet{ashwin2012tipping} identified
a few mathematical mechanisms behind the observed phenomena,
attempting a classification:
\begin{itemize}
\item Bifurcation-induced tipping (Slow passage through a bifurcation)
\item Noise-induced tipping (Transition between attractors due to
  random fluctuations)
\item Rate-induced tipping (Failure to track a continuously
  changing quasi-steady state).
\end{itemize}
This paper studies how a system that is close to a rate-induced
  tipping event behaves under the influence of additive noise. We look
  at a prototypical system, the saddle-node normal form with additive
  noise and a ramped shift of the equilibrium as proposed by
  \citet{ashwin2012tipping}. (A more general definition and further
properties of rate-induced tipping are given by
\citet{ashwin2015parameter}.)  Two early-warning indicators that are
commonly used for bifurcation-induced tipping with noise are an
increase of autocorrelation and an increase in variance in observed time series of system outputs \cite{scheffer2009early}. The
  most common argument, why a generic output time series of a system
  approaching bifurcation-induced tipping should show an increase in
  autocorrelation and variance assumes that the bifurcation of the
  deterministic part is a saddle-node bifurcation. Far away from the
bifurcation one can think of the state of the system as the
  position of an overdamped particle at the bottom of a slowly
softening potential well. Any small perturbation or
  disturbance will relax back to the equilibrium with a large decay
  rate \cite{scheffer2012anticipating}. As the bifurcation is approached,
the potential well will become shallower, that is, the decay rate
  will decrease \cite{kuehn2015multiple}. Thus, any disturbance or perturbation will have an
  increased and more long-lasting effect such that, in the presence of
  noise, autocorrelation and variance in observed time series
increases. 

In practice, the early-warning indicators are used on
  observational time series data of systems where quantitatively
  accurate models are unavailable, such as palaeoclimate temperature
  and CO$_2$ proxies \cite{ditlevsen2010tipping} and lake eutrophication \cite{wang2012flickering}. In the
  cited cases the early-warning indicators were not used for
  prediction (as they were about events in the past), but as evidence
  for (or against) the presence of underlying tipping mechanisms. For
  example, \citet{dakos2008slowing} used an increase in
  autocorrelation in a sequence of palaeoclimate time series as
  evidence for bifurcation-induced tipping while
  \citet{ditlevsen2010tipping} used the absence of the increase in
  variance (and the inconclusive behavior of the autocorrelation) as
  evidence that the Dansgaard-Oeschger events are a case of
  noise-induced tipping. 
Similarly, the presence of early-warning indicators in simulation data
from a global circulation model showing a collapse of the Atlantic
Meridional Overturning circulation (AMOC) through freshwater input was
used as evidence for a bifurcation-induced tipping event
\cite{boulton2014early}. 
All of the cited studies base their arguments on the knowledge
  that, close to a bifurcation- or noise-induced tipping point the
  system (after de-trending) behaves like an Ornstein-Uhlenbeck (OU)
  process. For the OU process one can infer from observed
  autocorrelation and variance the underlying linear decay rate, thus,
  permitting conclusions about the approach (or lack of it) of the
  equilibrium to a saddle-node bifurcation. This paper studies the effect of noise on the
third mechanism from the list by \citet{ashwin2012tipping},
rate-induced tipping, with the goal to aid identification of this type
from time series.

In contrast to bifurcation-induced tipping, rate-induced tipping is failure of the system to track the continuously changing quasi-steady state \cite{ashwin2012tipping}. Unlike bifurcation-induced tipping, at each moment in time there exists a stable (quasi-)equilibrium but the rate at which this steady state shifts determines whether the system tips or not.

The effect of rate-induced tipping has been described only
  relatively
  recently. 
In particular, within climate science
\citet{wieczorek2011excitability} considered a model for carbon
  storage and release in peatland soil, which showed the compost
bomb instability. In their model an increase in temperature above a
critical rate results in a release of carbon into the atmosphere from
combustion of compost heaps. A higher CO$_{2}$ concentration in the
atmosphere, creates further warming and thus triggering a positive
feedback loop within the system \cite{luke2011soil}. This is an
example of rate-induced tipping as for every fixed atmospheric
temperature there exists a globally stable steady state but the
rapidity of the temperature increase causes sharp peaks of carbon
  release. Other examples of rate-induced tipping include the
switching off of the AMOC due to the rate of increase of CO$_{2}$ in
the atmosphere \cite{stocker1997influence}. \citet{scheffer2008pulse}
find in a plant-herbivore model the critical rates of plant growth
causing a rate-induced transition from a herbivore controlled state to
a vegetated state.

Rate-induced tipping is not associated to a loss of stability of
equilibrium and thus cannot be explained using stability theory
for equilibria\cite{perryman2014adapting}. An appropriate
analogue to the``overdamped particle in a softening
well" illustration for bifurcation-induced tipping is to
  think of an overdamped particle in a moving well. 
In contrast to
bifurcation-induced tipping, the shape of the potential well remains
constant but instead shifts at varying rates. The faster the shift the
further the particle drifts away from the bottom of the well, up the
side and thus closer to the saddle and escaping. Hence, there is no
change in stability of the potential well, only the location of where
the state is in terms of the potential. As a consequence,
\citet{ashwin2012tipping} remarked there is no reason to assume why
the early-warning indicators: autocorrelation and variance can still
give useful predictions.

This paper builds on the work of \citet{ashwin2012tipping}, which
introduced a prototypical model for deterministic rate-induced
tipping. We will consider the effect of additive white noise on
  this prototype model of rate-induced tipping. This models
fluctuations/uncertainties that exist in various systems, for example
the climate system. It also permits us to study early-warning
  indicators. The aim of this paper is to demonstrate that
autocorrelation and variance will show an increase. However, this increase occurs with a delay, which is related to a delay in the actual tipping.

The paper is structured as follows: Section~\ref{sec:deterministic
  r-tipping} describes the basic properties of the
  deterministic prototype model for rate-induced tipping introduced by
  \citet{ashwin2012tipping}.  Section~\ref{sec:apparent failure of
  early-warning} explores the apparent delay of the early-warning
indicators for noise and rate-induced tipping. In
Section~\ref{sec:most likely time of tipping} we set up a
  boundary-value problem for most likely tipping paths, the sequence of continuation steps to solve this boundary-value problem are presented in Section~\ref{sec:steps}. In Section~\ref{sec:optimal} the most likely tipping path is discussed for a fixed set of system parameters, and in Section~\ref{sec:timing} analysis of most likely paths for all relevant system parameters using numerical continuation is covered. Section~\ref{sec:gendelay} discusses results of delay in the context of autonomous systems before,
Section~\ref{sec:conclusions} presents some concluding remarks.

\section{The deterministic backbone --- a prototype for rate-induced tipping}
\label{sec:deterministic r-tipping}

\begin{figure}[ht]
        \centering       
                \includegraphics[scale = 0.3]{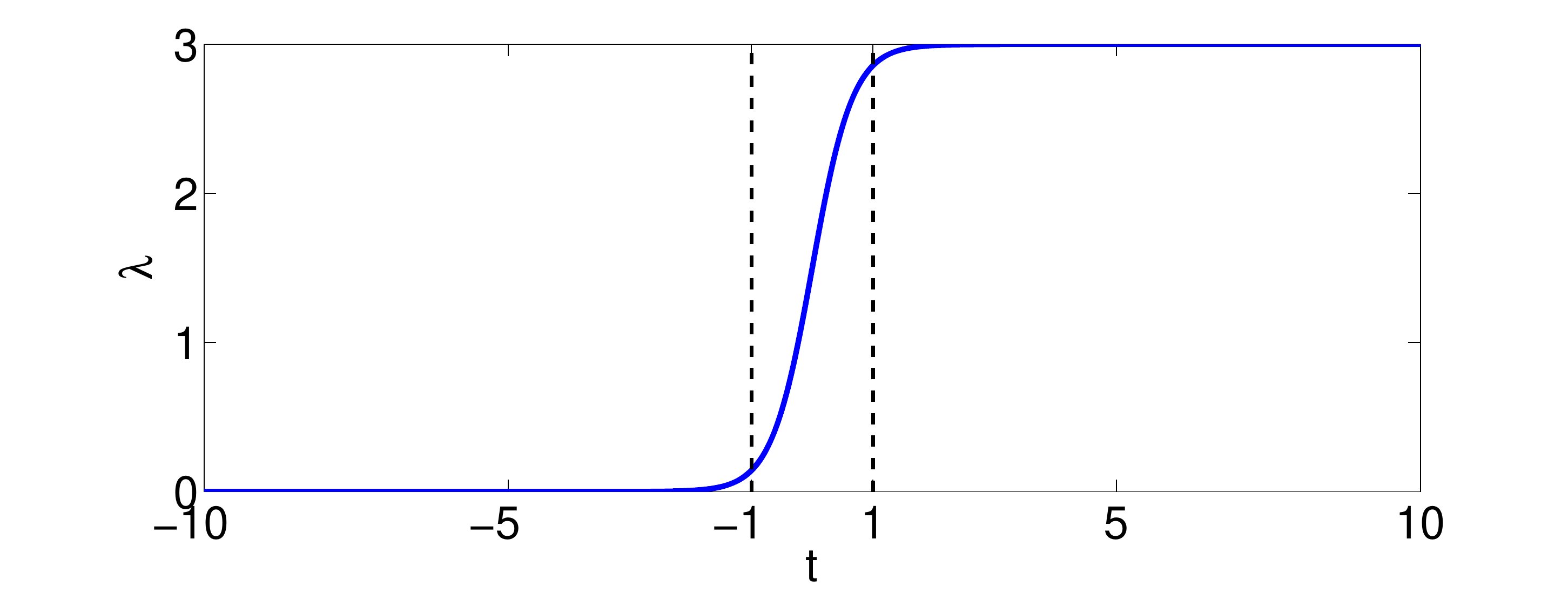}
                \caption{Time profile of shift parameter $\lambda(t)$, equation (\ref{lambda eq}), \change{where the black dashed lines indicate the transition period} for $\lambda_{\max} = 3$, $\epsilon = 1$}
                \label{ramp}
\end{figure}

A prototype model for rate-induced tipping was introduced by
\citet{ashwin2012tipping}. 
\change{The model is a scalar ordinary differential equation (ODE) for the variable $x(t)\in\R$:}
\begin{eqnarray}
\dot{x} = f(x,\lambda) = (x + \lambda)^{2} - 1\mbox{.} 
\label{rtip ODE initial}
\end{eqnarray}
\noindent \change{which, is the normal form for the saddle-node bifurcation.}
We have set the normal form parameter equal to $1$
w.l.o.g. (corresponding to a choice of scale for $x$ and time). The
ODE \eqref{rtip ODE initial} has two
$\lambda$-dependent families of equilibria, one stable at
$x_\mathrm{eq}^{(s)}(\lambda) = -\lambda - 1$ and one unstable at
$x_\mathrm{eq}^{(u)}(\lambda) = -\lambda + 1$. The equilibria are
separated by a distance of $2$.  These families of equilibria
$x_\mathrm{eq}^{(s)}(\lambda)$ and $x_\mathrm{eq}^{(u)}(\lambda)$ form
straight lines in the $(\lambda,x)$ - plane and will be referred to as
$W^{s}_{0}$ and $W^{u}_{0}$ respectively (see Figure \ref{Deterministic
  plots}). Equation~\eqref{rtip ODE initial} is the saddle-node
normal form shifted by $\lambda$, for which we assume dependence on
time in the form of a ramp (see Figure \ref{ramp}):
\begin{eqnarray}
\lambda(t) = \dfrac{\lambda_{\max}}{2}\bigg[\tanh\bigg(\dfrac{\lambda_{\max}\epsilon t}{2}\bigg) + 1\bigg]
\label{lambda eq}
\end{eqnarray} 
where $\lambda_{\max}$ (distance) and $\epsilon$ (speed) are the shape
parameters of the ramp-like shift. 
The time-derivative of $\lambda(t)$
is of most interest here as this determines the rate of shift for the
(now quasi-) equilibria $x_\mathrm{eq}^{(s)}$ and
$x_\mathrm{eq}^{(u)}$. The time derivative of $\lambda$ is
\begin{eqnarray}
\dfrac{\mathrm{d}\lambda}{\mathrm{d}t} = \dfrac{\epsilon\lambda_{\max}^{2}}{4}\bigg[\sech^{2}\bigg(\dfrac{\lambda_{\max}\epsilon t}{2}\bigg)\bigg]=\epsilon\lambda(\lambda_{\max}-\lambda)\mbox{.}
\label{lambda prime}
\end{eqnarray}
This time derivative reaches its maximum at $t=0$ and so, for a fixed
ramp height $\lambda_{\max}$, $\epsilon$ is directly proportional to
the maximal rate of shift at $t=0$. 


We note that \eqref{lambda prime} is also an ODE for $\lambda$ such that the prototype model can be considered as a two-dimensional ODE in the $(x,\lambda)$ phase plane (as done by \citet{ashwin2012tipping}):
\begin{eqnarray}
\label{rtip ODE}
\dot{x} &&= f(x,\lambda(t)) = (x + \lambda(t))^{2} - 1 \\
\label{rtip lambda dot init}
\dot{\lambda} &&= \epsilon\lambda(\lambda_{max} - \lambda)
\end{eqnarray}  
Notice that \eqref{rtip lambda dot init} is coupled to \eqref{rtip
  ODE}, but there is no coupling in the other direction.
\begin{figure}[ht]
        \centering
        \subcaptionbox{\label{TP eps small}}[0.45\linewidth]
                {\includegraphics[scale = 0.3]{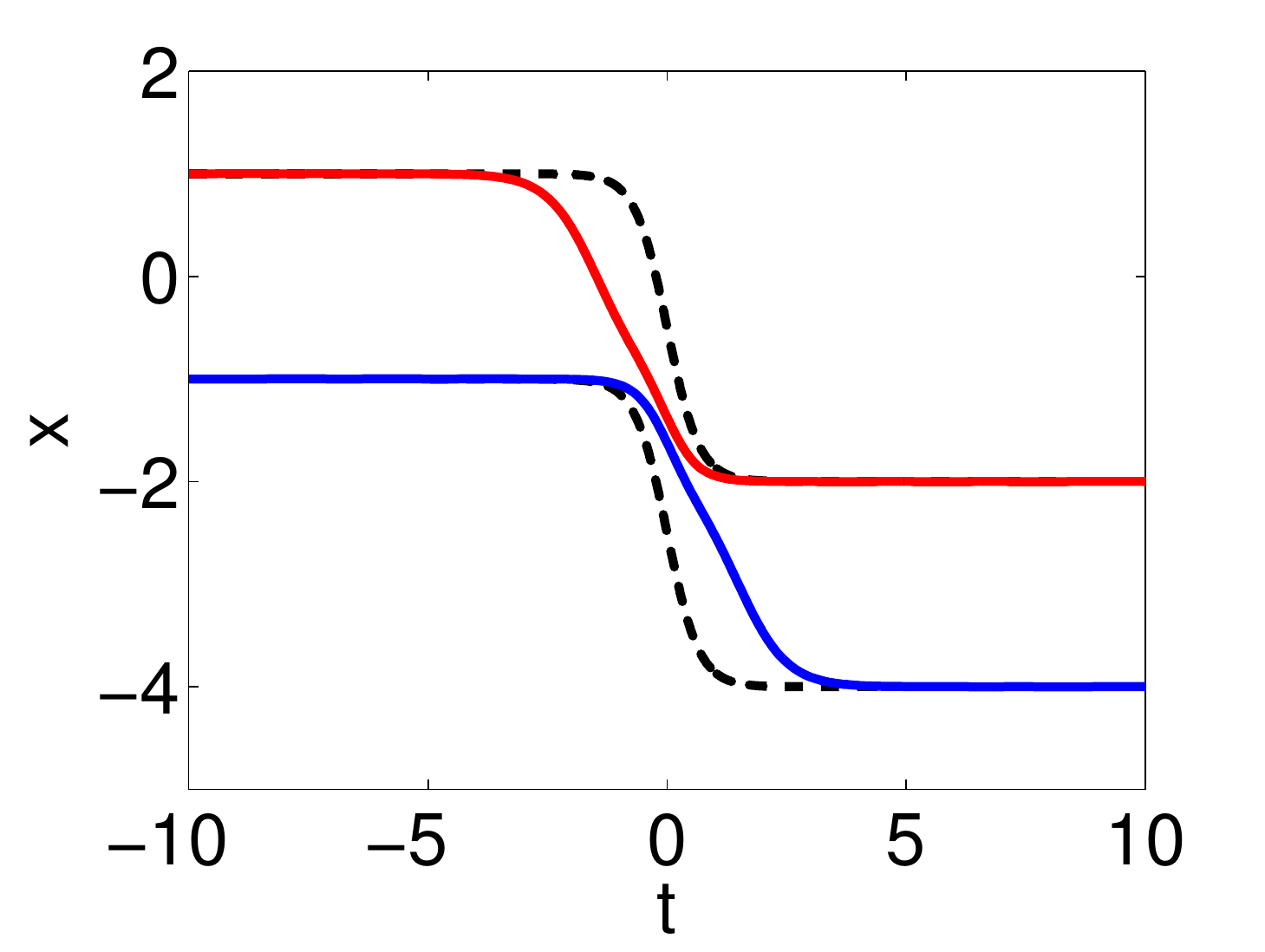}}
       \hfill 
        \subcaptionbox{\label{PP eps small}}[0.45\linewidth]
                {\includegraphics[scale = 0.3]{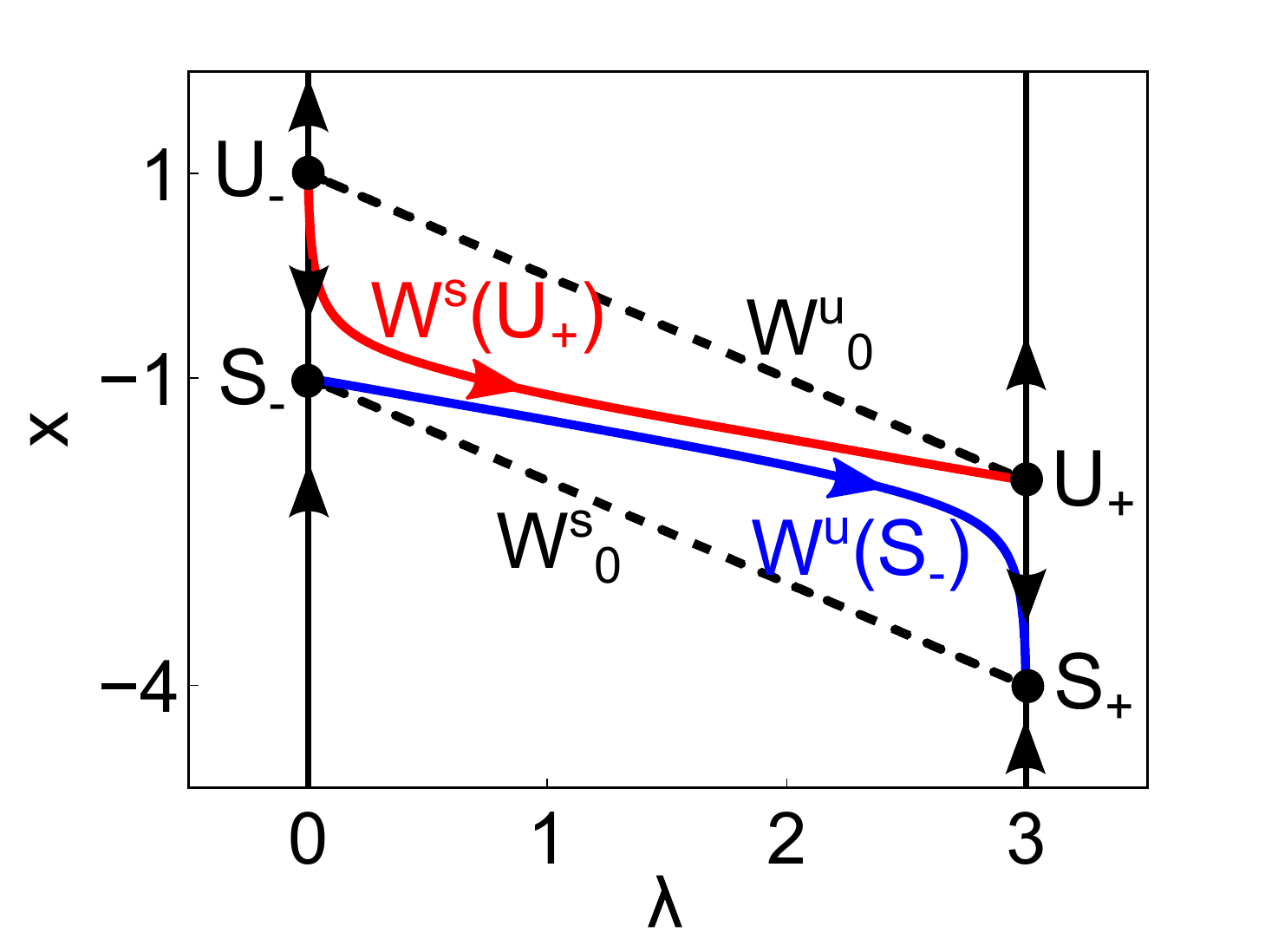}}
        ~ 
          \\
          \subcaptionbox{\label{TP heteroclinic}}[0.45\linewidth]
                {\includegraphics[scale = 0.3]{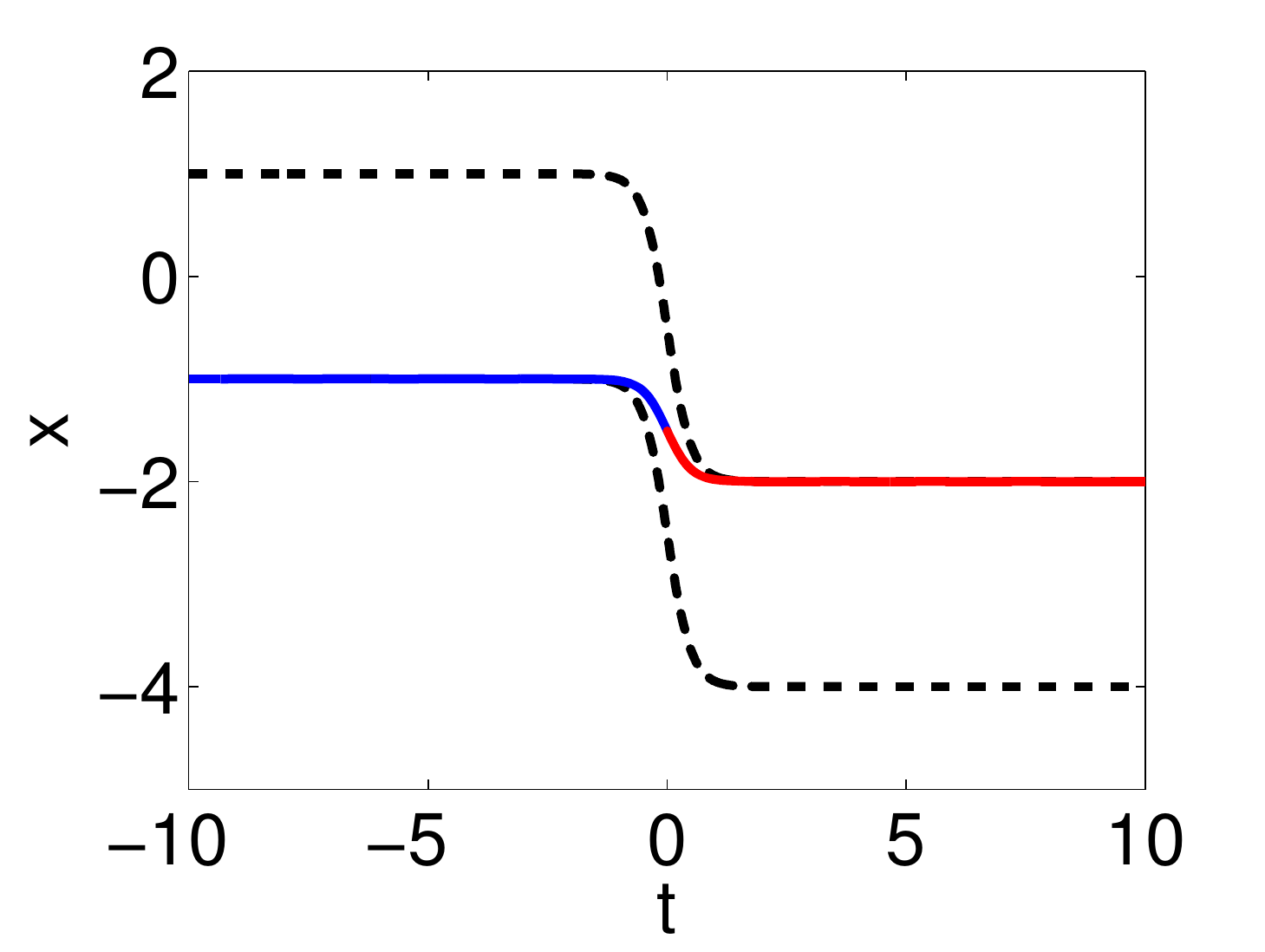}}
        \hfill 
        \subcaptionbox{\label{PP heteroclinic}}[0.45\linewidth]
                {\includegraphics[scale = 0.3]{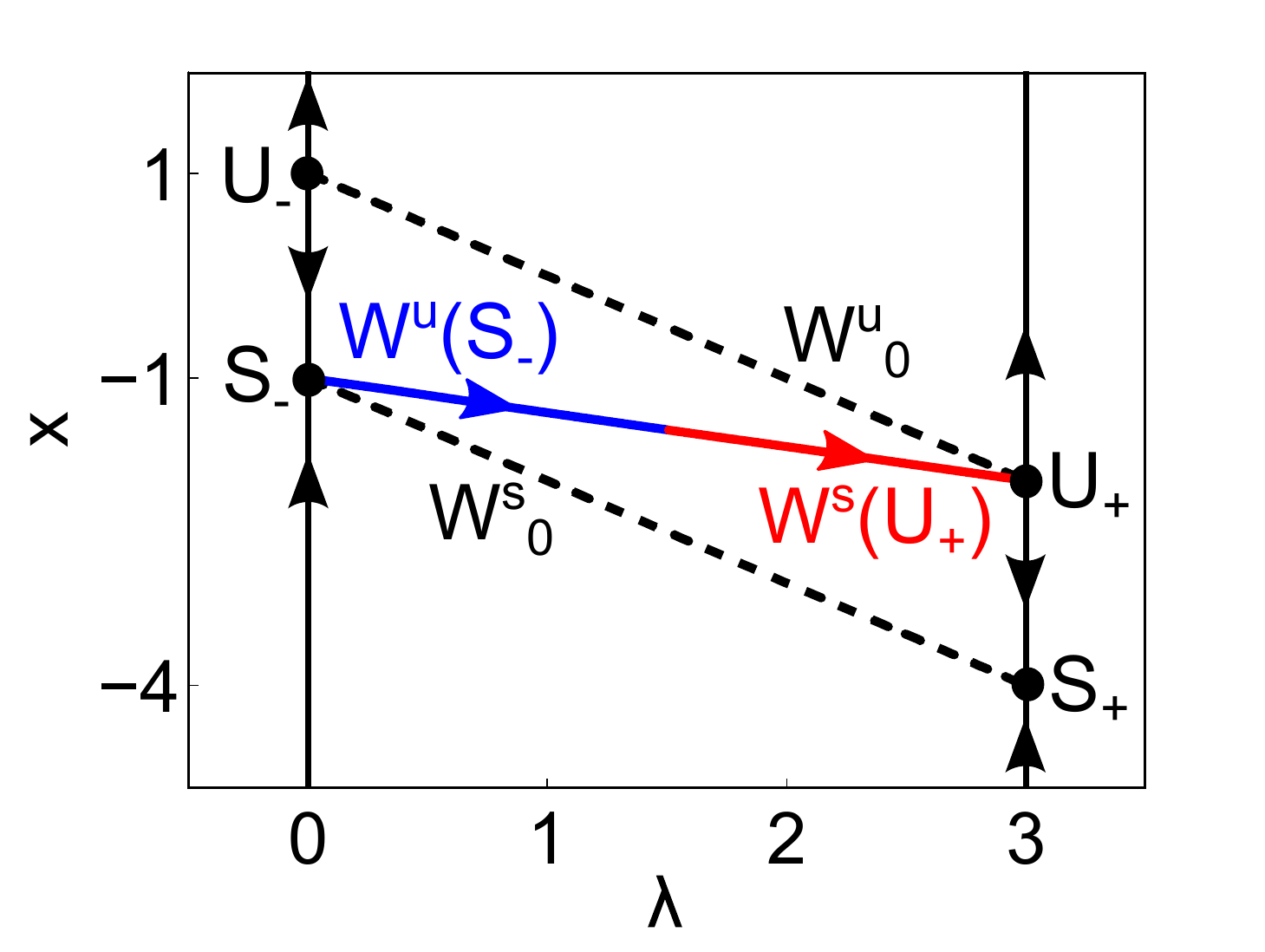}}
        ~ 
          \\
          \subcaptionbox{\label{TP eps big}}[0.45\linewidth]
                {\includegraphics[scale = 0.3]{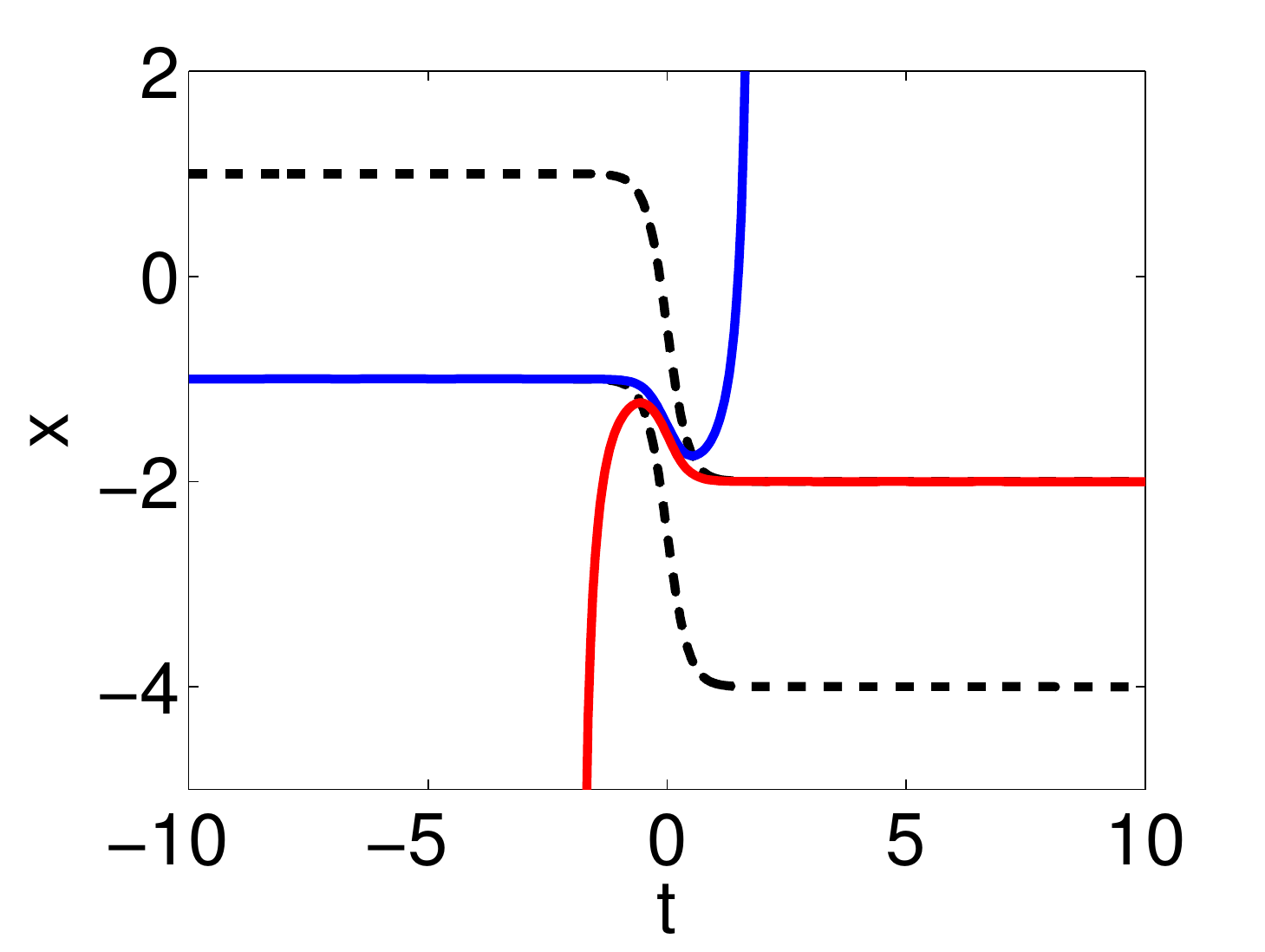}}
       \hfill 
        \subcaptionbox{\label{PP eps big}}[0.45\linewidth]
                {\includegraphics[scale = 0.3]{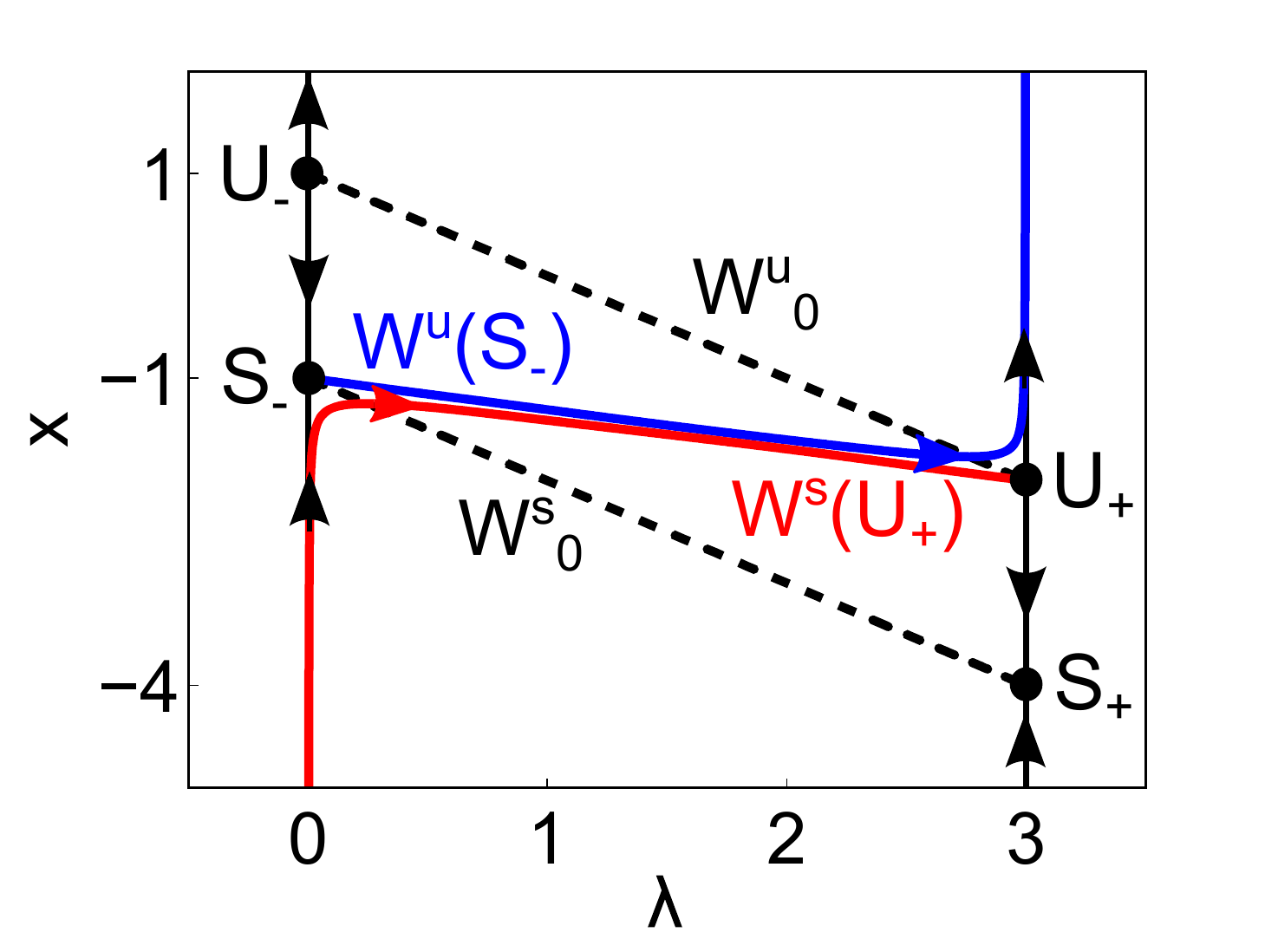}}
        ~ 
        \caption{Time profiles (a),(c),(e) and phase planes
          (b),(d),(f) of system \eqref{rtip ODE}--\eqref{rtip lambda dot
            init} for $\epsilon < \epsilon_{c}$ - (a),(b), $\epsilon
          = \epsilon_{c}$ - (c),(d) and $\epsilon > \epsilon_{c}$ -
          (e),(f). Black dashed curves are the stable $W^{s}_{0}$ and
          unstable $W^{u}_{0}$ equilibria in the limit $\epsilon = 0$,
          blue and red curves are the unstable and stable manifolds,
          $W^{u}(S_{-})$ and $W^{s}(U_{+})$,
          respectively ($\lambda_{\max}=3$).}\label{Deterministic plots}
\end{figure}
Figure \ref{Deterministic plots} displays all qualitatively different
phase portraits possible for \eqref{rtip ODE}--\eqref{rtip lambda dot
  init} in panels \ref{PP eps small}, \ref{PP heteroclinic} and
\ref{PP eps big}. The system has $4$ equilibria, $S_-$ (a saddle) and
$U_-$ (a source) on the $\lambda=0$ line, and $S_+$ (a sink) and $U_+$
(a saddle) on the $\lambda=\lambda_{\max}$ line. The upper and lower
black dashed lines represent the family of unstable $W^{u}_{0}$ and
stable $W^{s}_{0}$ quasi-equilibria in the limit $\epsilon = 0$,
respectively. The blue curve is the unstable manifold $W^{u}(S_{-})$
of the saddle $S_{-}$, and the red curve is the stable manifold
$W^{s}(U_{+})$ of the saddle $U_{+}$. The panels \ref{TP eps small},
\ref{TP heteroclinic} and \ref{TP eps big} show the time profiles for
$x$ on the invariant manifolds $W^{u}(S_{-})$ and $W^{s}(U_{+})$
(using the same color coding).

One can see that the time profile and the phase portrait of the unstable
manifold $W^{u}(S_{-})$ and the stable manifold $W^{s}(U_{+})$
(Figures \ref{TP eps small} and \ref{PP eps small}) deviate
increasingly from the quasi-equilibrium families $W^{u}_{0}$ and  $W^{s}_{0}$ for
increasing $\epsilon$.
For small $\epsilon$, $W^{u}(S_{-})$ is close to $W^{s}_{0}$, but, for
increasing $\epsilon$, $W^{u}(S_{-})$ moves further apart from
$W^{s}_{0}$. The unstable manifold $W^{u}(S_{-})$ converges for
$t\to\infty$ to the stable node $S_{+}$ for
$\epsilon<\epsilon_c$. The red curve is the stable manifold
$W^{s}(U_{+})$, which 
forms a separatrix partitioning the plane into two regions. In the
region below the separatrix all trajectories are attracted to the
stable node $S_{+}$, but in the region above the repelling stable
manifold $W^{s}(U_{+})$, all trajectories will escape to $+\infty$ in
finite time. Notice that the two manifolds $W^{u}(S_{-})$ and
$W^{s}(U_{+})$ are closest at $\lambda=\lambda_{\max}/2$, when the
time-derivative $\lambda$ (equation \eqref{lambda prime}) is at its
maximum. This is due to the reflection symmetry within the system
\eqref{rtip ODE}, \eqref{rtip lambda dot init}

\begin{eqnarray*}
\begin{bmatrix}
x-x_{c} \\
\lambda-\lambda_{c}
\end{bmatrix}
\rightarrow
\begin{bmatrix}
x_{c}-x \\
\lambda_{c}-\lambda
\end{bmatrix}
\end{eqnarray*}
around the point $(x_{c},\lambda_{c}) = (-1.5,1.5)$.

At a critical $\epsilon$, denoted $\epsilon_{c}$, $W^s(U_+)$ and
$W^u(S_-)$ form a heteroclinic connection between the two saddles
$S_{-}$ and $U_{+}$, as depicted in Figures \ref{TP heteroclinic},
\ref{PP heteroclinic}.
\citet{Perryman2015tipping} observed that the critical value $\epsilon_c$ equals $4/3$ and that the connecting orbit is the line
\begin{eqnarray}
x = -\dfrac{\lambda}{3} - 1
\label{heteroclinic}
\end{eqnarray} 
in the phase plane.  For $\epsilon>\epsilon_{c}$, $W^{u}(S_{-})$ and
$W^{s}(U_{+})$ change their arrangement, as displayed by Figures
\ref{TP eps big}, \ref{PP eps big}. The unstable manifold
$W^{u}(S_{-})$ no longer converges to the stable node $S_{+}$ such
that trajectories from all initial conditions close to $S_-$ with
$\lambda>0$ diverge ($x(t)\to+\infty$).  In this case the parameter
$\lambda$ is shifted at a rate that is too large for the unstable
manifold $W^{u}(S_{-})$ to track the quasi-steady state $W^s_0$. For
$\epsilon>\epsilon_c$ but close to $\epsilon_c$ this escape does not
occur until $\lambda$ is close to $3$ such that one would observe the escape only
when $\lambda$ is coming to rest again, and so there appears to be a lag in
the timing of escape.

In summary, for this prototype model rate-induced tipping corresponds
to a global bifurcation at parameter value $\epsilon_{c}$, a
heteroclinic connection from the stable equilibrium before the
parameter ramp to the unstable equilibrium after the ramp.


\section{Delay of early-warning indicators and delay of tipping}
\label{sec:apparent failure of early-warning}

For the remainder of this paper we will consider the following
scenario: the speed of the parameter ramp $\epsilon$ is less than its
critical value $\epsilon_{c}=4/3$ such that without noise the system
will not tip. The influence of noise, which we add to the dynamics
\eqref{rtip ODE} of $x$, will cause the system to tip with a certain
probability. We can control this probability by varying noise
intensity and $\epsilon$. We choose our parameters such that an escape
of $x$ from $W^s_0$ beyond $W^u_0$ to $+\infty$ is extremely unlikely
for $t$ far away from $0$ (and, thus, $\lambda$ far away from
$\lambda_{\max}/2$). We expect this escape probability to increase
during the ramp of $\lambda$ (for $t\approx 0$).

The realizations of $x$ for the prototype system \eqref{rtip ODE}--\eqref{rtip lambda dot init} are governed by the stochastic differential equation (SDE):
\begin{eqnarray}
\mathrm{d}X_{t} = [(X_{t} + \lambda(t))^2 - 1]\mathrm{d}t + \sqrt{2D}\mathrm{d}W_{t}
\label{SDE}
\end{eqnarray}
where $W_{t}$ is standard Brownian motion. The intensity of the noise is given by $\sqrt{2D}$ where $D$ is the diffusion coefficient. The probability density $P(x,t)$ of the random variable $X_{t}$ in the SDE \eqref{SDE} is governed by the Fokker-Planck equation; a linear partial differential equation (PDE):
\begin{eqnarray}
\dfrac{\partial P(x,t)}{\partial t} = D\dfrac{\partial^{2}P(x,t)}{\partial x^{2}} - \dfrac{\partial}{\partial x}\bigg(f(x,t)P(x,t)\bigg)\mbox{,}
\label{FPE}
\end{eqnarray}
which includes the diffusion coefficient $D$ and drift term $f(x,t)=(x+\lambda(t))^2-1$. \change{Applying Dirichlet boundary conditions at some $[x_{\mathrm{start}},x_{\mathrm{end}}]$ will cause the probability density to decay over time as realizations escape the domain. The probability of escape, $p_{\mathrm{esc}}(t_{n})$ at time step $t_{n}$ is therefore defined as:

\begin{equation*}
p_{\mathrm{esc}}(t_{n}) = 1 - \frac{\int P(x,t_{n})\mathrm{d}x}{\int P(x,t_{n-1})\mathrm{d}x}
\end{equation*}

} 

\change{In addition,} we use \eqref{FPE} to compute two characteristic quantities of the
density $P(x,t)$ for \eqref{SDE}, the lag-1 autocorrelation and the
variance, shown in Figure~\ref{Early-warning indicators}. These
quantities are commonly monitored in time series where one suspects an
underlying parameter drift that approaches a bifurcation-induced
(specifically saddle-node induced) tipping point. Both,
autocorrelation and variance, should increase along the time series as
the parameter comes closer to its saddle-node value (see
\citet{williamson2015detection} for other cases such as Hopf
bifurcation). But what happens for the rate-induced tipping model?

\begin{figure}[ht]
        \centering
        \subcaptionbox{Autocorrelation\label{Autocorrelation}}[0.45\linewidth]
                {\includegraphics[scale = 0.3]{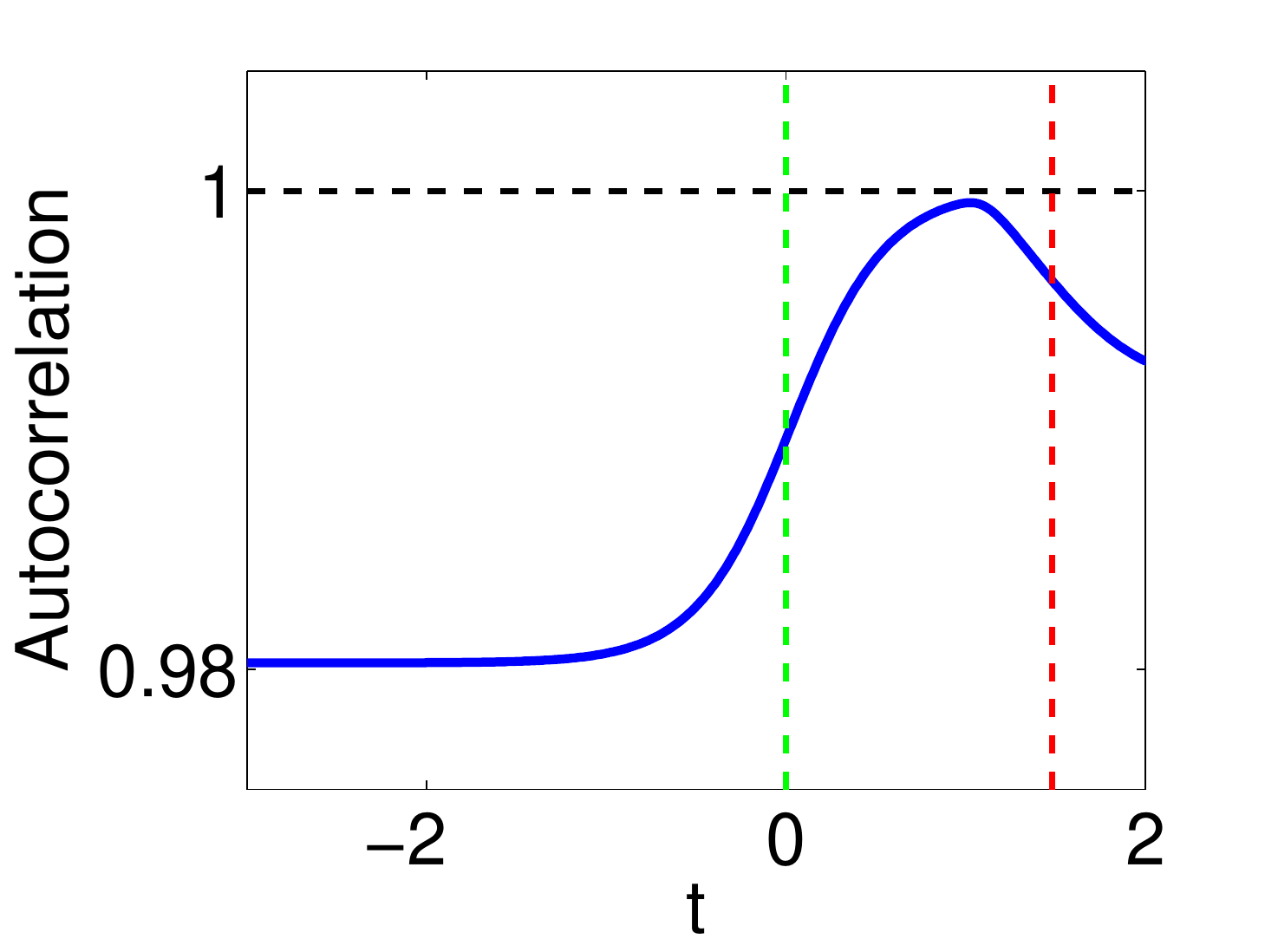}}
        \hfill 
        \subcaptionbox{Variance\label{Variance}}[0.45\linewidth]
                {\includegraphics[scale = 0.3]{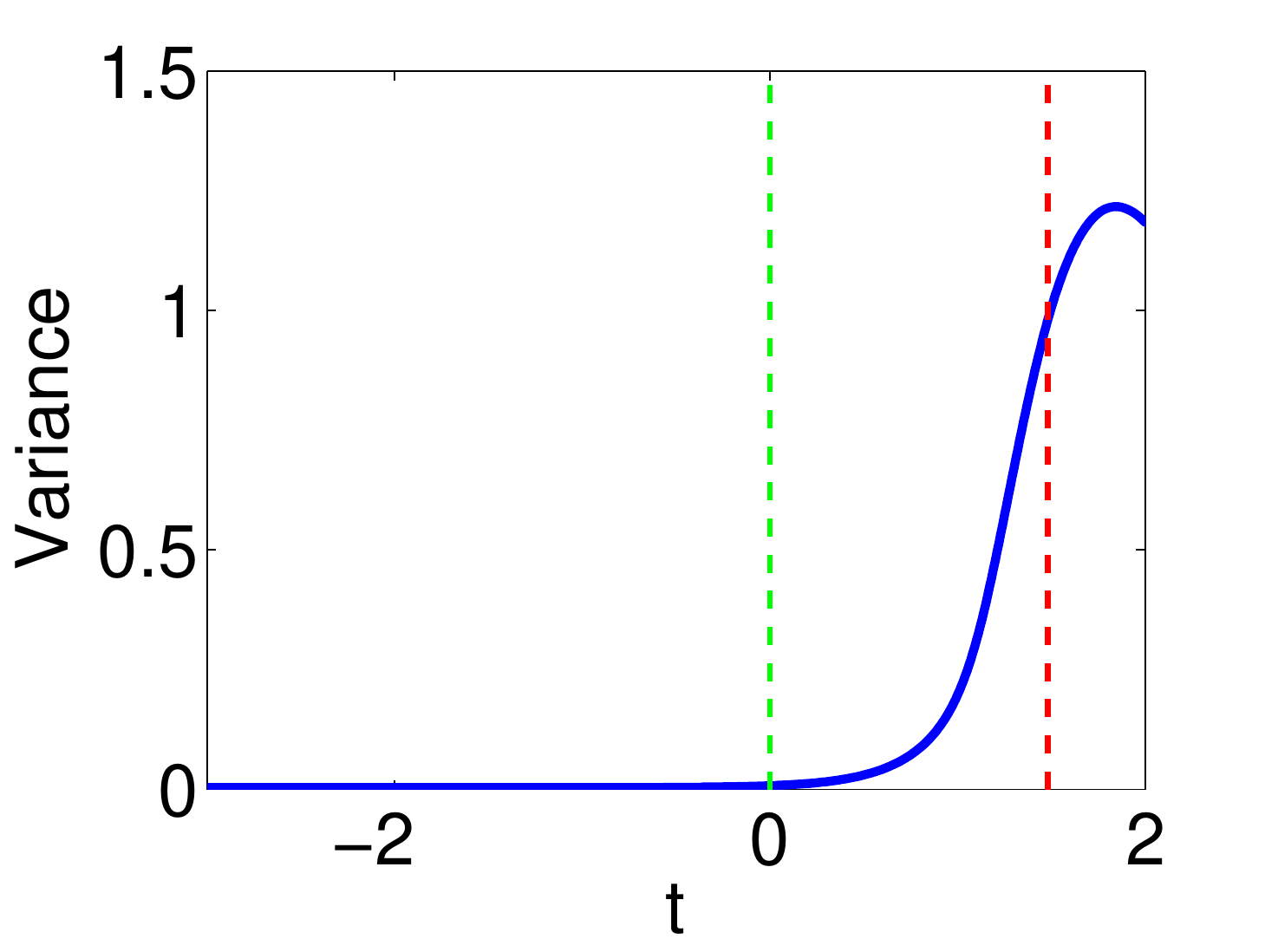}}
        ~ 
        \caption{Traditional early-warning indicators; lag-1
          autocorrelation and variance show a delayed response for
          rate-induced tipping if tipping point is assumed to be at the time $t=0$ (green dashed line), the closest encounter of the stable
  and unstable manifolds, $W^{u}(S_{-})$ and $W^{s}(U_{+})$. However, both indicators increase before the tipping point if assumed to be at $t \approx 1.5$ (red dashed line) calculated from the peak of the escape rate in Figure \ref{Escape rate}. Parameters: $\epsilon = 1.25$ and
          $D = 0.008$ ($\Delta t=0.01$
          for panel \ref{Autocorrelation})}\label{Early-warning
              indicators}
\end{figure}

The lag-1 autocorrelation $a_{n}$ is defined to be the correlation
between successive $X_{n\Delta t}$, separated by a time
step $\Delta t$ (we choose $\Delta t= 0.01$):

\begin{eqnarray}
a_{n} = \dfrac{\mathrm{Cov}(X_{(n-1)\Delta t},X_{n\Delta t})}{
  \sqrt{\mathrm{Var}(X_{(n-1)\Delta t})\mathrm{Var}(X_{n\Delta t})}}
\label{autocorrelation}
\end{eqnarray}
where $X_{n\Delta t}$ is the solution of \eqref{SDE} with density
$P(\cdot,n\Delta t)$ at time step $n\Delta t$.

The initial condition for \eqref{FPE} is the stationary density of
\eqref{SDE} with \change{$\lambda_{0} = \lambda(t_{0})$} restricted to the \change{fixed} domain $x\in[x_{\mathrm{start}},x_{\mathrm{end}}]=[-6,2]$, which
corresponds to the assumption that the ramp-up of $\lambda$ starts
from a stationary state. \change{See Appendix A(a) for a study of dependence on $x_{\mathrm{end}}$. For this stationary starting point the system can be approximately modeled by the Ornstein-Uhlenbeck process

\begin{equation*}
\mathrm{d}X_{t} = -\theta X_{t}\mathrm{d}t + \sqrt{2D}\mathrm{d}W_{t}
\end{equation*}

\noindent where $\theta = -f'(-1,0) = 2$ is the decay rate at $S_{-}$. The Ornstein-Uhlenbeck process has autocorrelation and variance given by \cite{aalen2008survival}:}

\begin{eqnarray*}
\text{Autocorrelation:} \qquad &&a = \exp(-\theta\Delta t) \approx 1 - \theta\Delta t \\
\text{Variance:} \qquad &&V = \dfrac{D}{\theta}
\end{eqnarray*}
where we set
$\Delta t = 0.01$, thus, giving $a=0.98$ and $V=0.004$ in Figure~\ref{Early-warning
  indicators}.
  
\change{We highlight that for starting at $t_{0} = -\infty$ the system will tip with probability one before the ramping shift begins. The time of tipping for a stationary system is approximated by Kramers' time, $\tau_{K}$ \cite{berglund2004noise}, $(U(x) = -\int f(x)\mathrm{d}x)$: 

\begin{equation*}
\tau_{K} = C\exp\bigg(\dfrac{\Delta U}{D}\bigg)
\end{equation*}

\noindent where $\Delta U$ is the height of the potential barrier and
the prefactor $C$ depends on the curvature at the minimum and maximum
of the potential. In our case for equation \eqref{SDE}, we have constant values for $\Delta U = 4/3$ and $C = \pi$. \change{We consider the regime where the probability of
  escape from the well, $p_{\mathrm{esc}}(t)$, increases by an order of
  magnitude during the ramping of the parameter $\lambda$:
  \begin{displaymath}
    \frac{1}{\tau_K}\ll \max_{t\in\R}p_\mathrm{esc}(t)\mbox{.}
  \end{displaymath}
  In this regime we expect the time $t$ for which the escape rate
  $p_\mathrm{esc}(t)$ is maximal} to occur at the time $t = 0$. This
corresponds to the closest encounter of the stable and unstable
manifolds, $W^{u}(S_{-})$ and $W^{s}(U_{+})$, due to the symmetry in
the deterministic part.}

Figure~\ref{Early-warning indicators} displays the lag-1 autocorrelation and variance for the time interval of most interest, namely $t\in[-3,2]$ when system \eqref{rtip ODE}, \eqref{rtip lambda dot init} is non-stationary. 
We observe that there is a
delay in the warning for approaching the tipping point, if we take
  the tipping point as the time $t=0$ (green dashed line). The autocorrelation (Figure
\ref{Autocorrelation}) has only just started to increase at $t =
0$. The variance (Figure \ref{Variance}) shows an even longer delay in the signal. It increases noticeably only after $t=0$. \citet{ditlevsen2010tipping} concluded for saddle-node
induced tipping that only the presence of both indicators, increase of
autocorrelation and variance, is sufficient evidence for the approach
of a tipping point. When applied to rate-induced tipping, one would
conclude initially that the warning will be significantly delayed from
when we would expect the tipping.  This warrants a systematic investigation to
see when escape is most likely in close encounters with rate-induced tipping.

\subsection{Escape rate over time} To investigate when the escape is likely to occur we initially consider the escape rate per unit time calculated using the
Fokker-Planck equation \eqref{FPE}. The escape rate per unit time is
defined as the fraction of realizations that cross a known threshold
curve $\tilde{x}(t)$ divided by the time step $\Delta t$. We choose a
threshold curve $\tilde{x}(t)$ (bright blue in Figure~\ref{Escape
  rate}) beyond which we classify a realization as having escaped. The
threshold curve is chosen such as $\tilde{x}(t) = x^{u}(t) + y$, where
$x^{u}(t)$ is the unique trajectory of the deterministic part of
\eqref{SDE} that starts at $x(-10) = x_{0} = -1$ (thus,
$(x^u(t),\lambda(t))$ is close to $W^u(S_-)$) and $y = 1.5$ is a fixed
sufficiently large deviation from
$x^u(t)$. Appendix~\ref{subapp:threshold} studies systematically how the
choice of threshold $\tilde x(t)$ affects the results.

\begin{figure}[ht]
        \centering
        \subcaptionbox{Time profile\label{TP escape}}[0.45\linewidth]
                {\includegraphics[scale = 0.3]{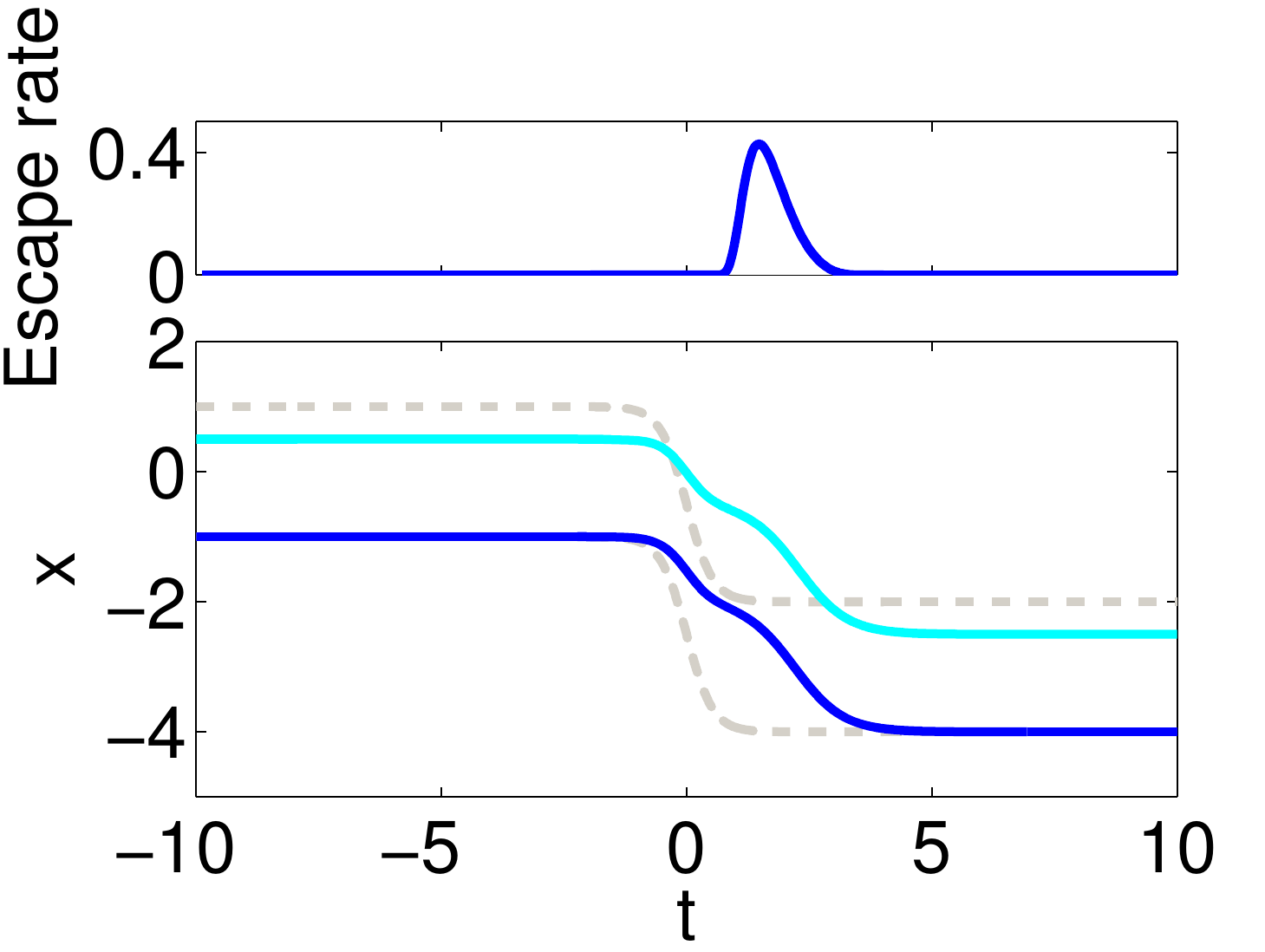}}
        \hfill 
        \subcaptionbox{Phase plane\label{PP escape}}[0.45\linewidth]
                {\includegraphics[scale = 0.3]{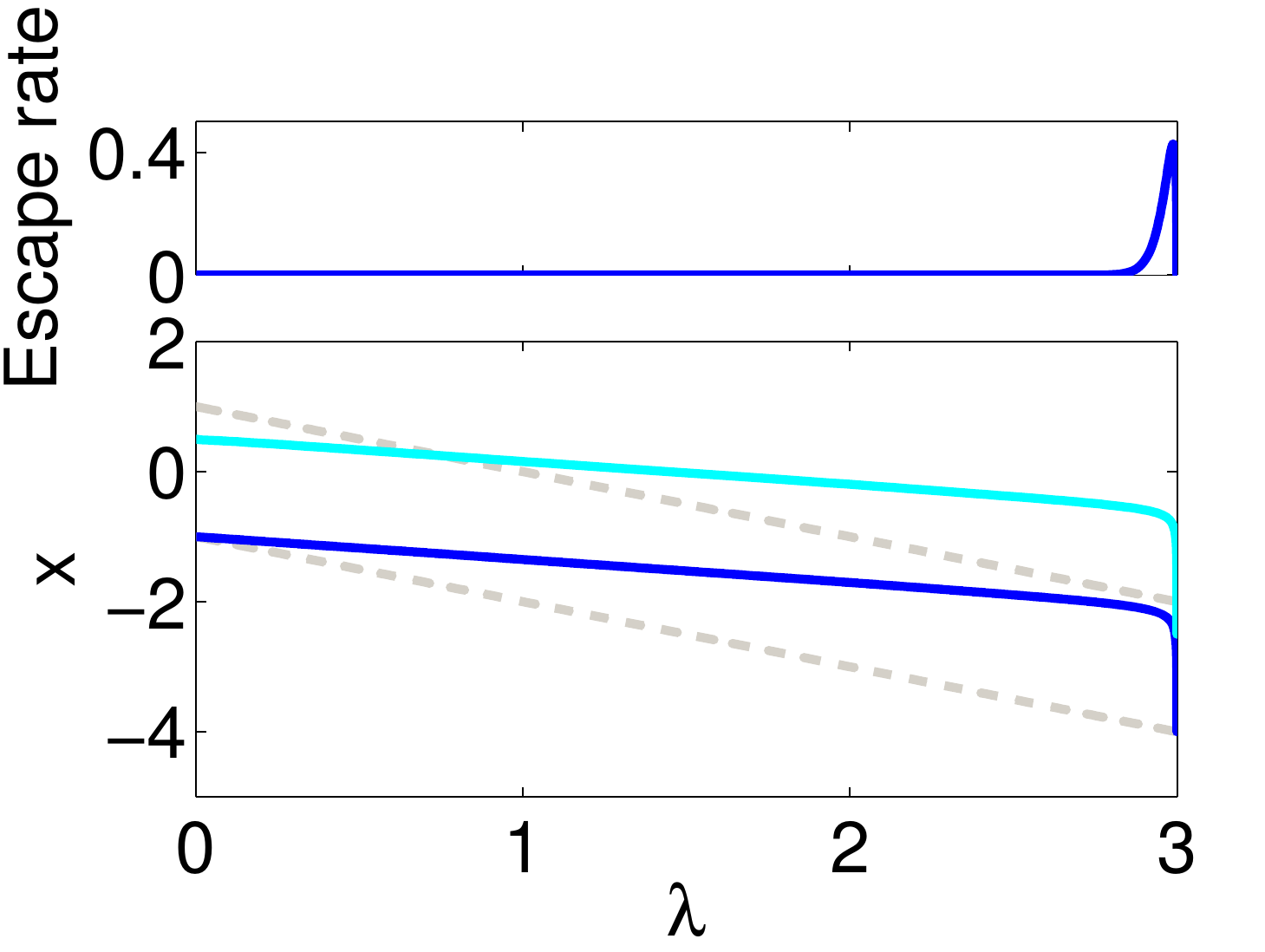}}
        ~ 
        \caption{Time profile and phase plane and escape rate obtained via \eqref{FPE} for $\epsilon = 1.25$ and $D = 0.008$. Dashed curves represent stable $W^{s}_{0}$ and unstable $W^{u}_{0}$ equilibria, the dark blue curve is the unstable manifold $W^{u}(S_{-})$, the bright blue curve is the threshold curve $\tilde{x}(t)$ for $y = 1.5$.}\label{Escape rate}
\end{figure}

Figure \ref{Escape rate} displays the time profile (\ref{TP escape})
and the phase portrait (\ref{PP escape}) of the deterministic
trajectory $x^u(t)$, the threshold curve $\tilde x(t)$ and the escape rate
(over time and versus $\lambda$) obtained via \eqref{FPE}.

In this example, we have chosen $\epsilon = 1.25$, which is close to
$\epsilon_{c}=4/3$, and a small noise level $D=0.008$. For this choice
of ramping speed parameter $\epsilon$, tipping would not occur in the
deterministic case ($D = 0$). However, with a noise level of $D =
0.008$, roughly \change{$36\%$} of realizations that start with initial
condition $x(-10) = x_{0}$ go on to escape. As in Figure
\ref{Deterministic plots}, the dark blue curve is the unstable
manifold $W^{u}(S_{-})$ (the deterministic solution $x^u(t)$, starting
at $x(-10) = x_{0}$, is extremely close to it). The bright blue curve
is our threshold curve $\tilde{x}(t)$ (see Figure \ref{Escape
  rate color} in Appendix~\ref{subapp:threshold} how the escape time depends on the threshold).

Figure~\ref{Escape rate} shows that the escape is most likely to occur
at about $t = 1.5$, hence, it is delayed, too. The red dashed line in Figure~\ref{Early-warning indicators} represents the most likely time of tipping, given by the peak of the escape rate in Figure~\ref{Escape rate}. Therefore taking this as our tipping point we see that there is an increase in both the autocorrelation and the variance on the approach to tipping. We conclude that for this example at least that the tipping is delayed and thus the early-warning signals can still offer forewarning of rate-induced tipping.  

For a systematic study of
how the time of most likely escape depends on the system parameters we
formulate a variational optimization problem for the optimal path of
escape.

\section{Most likely (optimal) escape paths --- The general variational problem}
\label{sec:most likely time of tipping}

In this section we will formulate the ODE boundary-value problem (BVP) determining locally most likely paths for escape. 

We define the most likely escape path as the path going from given
$x_0$ to a given $x_T$ in a time interval $[t_0,T_{\mathrm{end}}]$ maximizing
the functional
\begin{eqnarray}
F = \exp&&\bigg[\dfrac{U_{0} - U_{T}}{2D} - \int_{t_{0}}^{T_{\mathrm{end}}}\bigg(\dfrac{\dot{x}^{2}}{4D} + V_{s}\bigg)\mathrm{d}t\bigg]
\label{Functional}
\end{eqnarray}
along the path. The terms in $F$ are
\begin{align}
\nonumber
U(x,t) &= -\int f(\bar{x},t)\mathrm{d}\bar{x}\mbox{,}\\
\label{Vs equation}
V_{s}(x,t) &= \dfrac{1}{4D}\bigg(\dfrac{\partial U}{\partial x}\bigg)^{2} - \dfrac{1}{2}\dfrac{\partial^{2}U}{\partial x^{2}} - \dfrac{1}{2D}\dfrac{\partial U}{\partial t}\mbox{,}\\
\nonumber
  U_{0} &= U(x_{0},t_{0})\mbox{,} \quad U_{T} = U(x_{T},T_{\mathrm{end}})\mbox{.}
\end{align}

The quantity $U$ is the potential of the deterministic part $f$ of the
SDE \eqref{SDE} such that \eqref{SDE} can be written in terms
of $U(x,t)$:
\begin{eqnarray}
\mathrm{d}X_{t} = -\dfrac{\partial U(X_{t},t)}{\partial X_{t}}\mathrm{d}t + \sqrt{2D}\mathrm{d}W_{t}\mbox{.}
\label{SDE potential}
\end{eqnarray}
For a differentiable path $x$ the functional $F$ equals the
probability of a realization $X_t$ following a sequence of infinitesimally
small intervals $[x(k\Delta t)-\delta/2,x(k\Delta t)+\delta/2]$ in the
limit $0<\delta\ll\Delta t\ll1$ (up to a constant factor independent of $x$). Recall,
that the random variable $X_{t}$ had a probability density function
$P(x,t)$ given by the linear Fokker-Planck equation \eqref{FPE}, from which the functional $F$ is derived. \change{See \citet{zhang2008theory}, (p. 25-31) for a detailed derivation of the functional $F$ from equations \eqref{Functional}-\eqref{SDE potential} for a time independent potential $U(x)$ (a simple extension is made for a time dependent potential $U(x,t)$ in \citet{lin2011similarity} and \citet{ho2008perturbative}).} 


Assuming a fixed time interval $[t_0,T_{\mathrm{end}}]$ and fixed start and end
points $x_0$ and $x_T$, local critical points of $F$ are given by the
Euler-Lagrange equation, a 2nd order BVP
\cite{zhang2008theory}:
\begin{eqnarray}
\ddot{x} = 2D\dfrac{\partial V_{s}}{\partial x}(x,t)\mbox{,} 
\qquad \begin{cases}x(t_{0}) &= x_{0}\mbox{,} \\
x(T_{\mathrm{end}}) &= x_{T}\mbox{.} \end{cases}
\label{2nd order ODE}
\end{eqnarray} 

\begin{figure}[ht]
        \centering
        \subcaptionbox{\label{Tend continuation optimal path}}[0.45\linewidth]
                {\includegraphics[scale = 0.3]{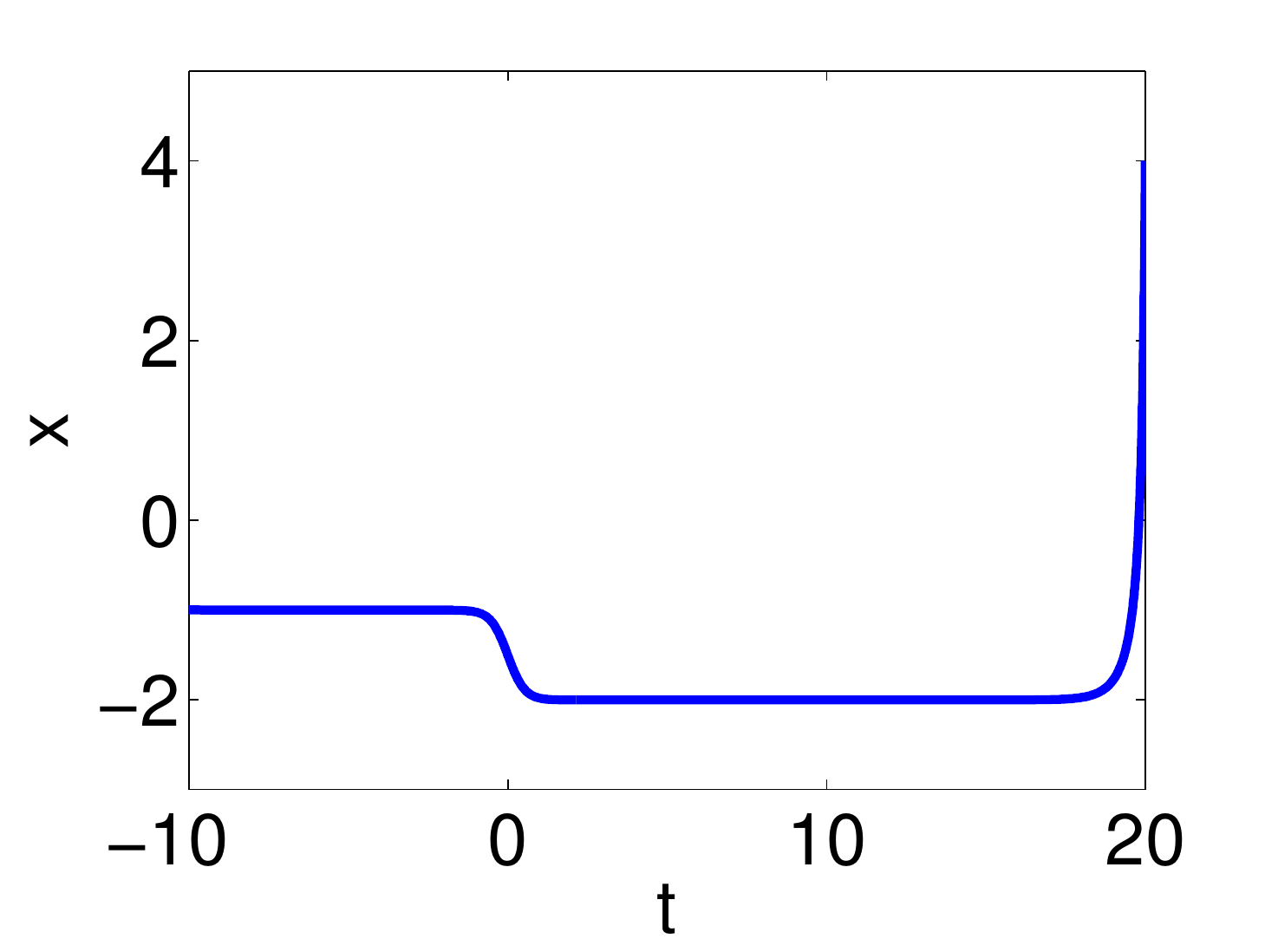}}
                \hfill
        \subcaptionbox{\label{Tend continuation 1}}[0.45\linewidth]
                {\includegraphics[scale = 0.3]{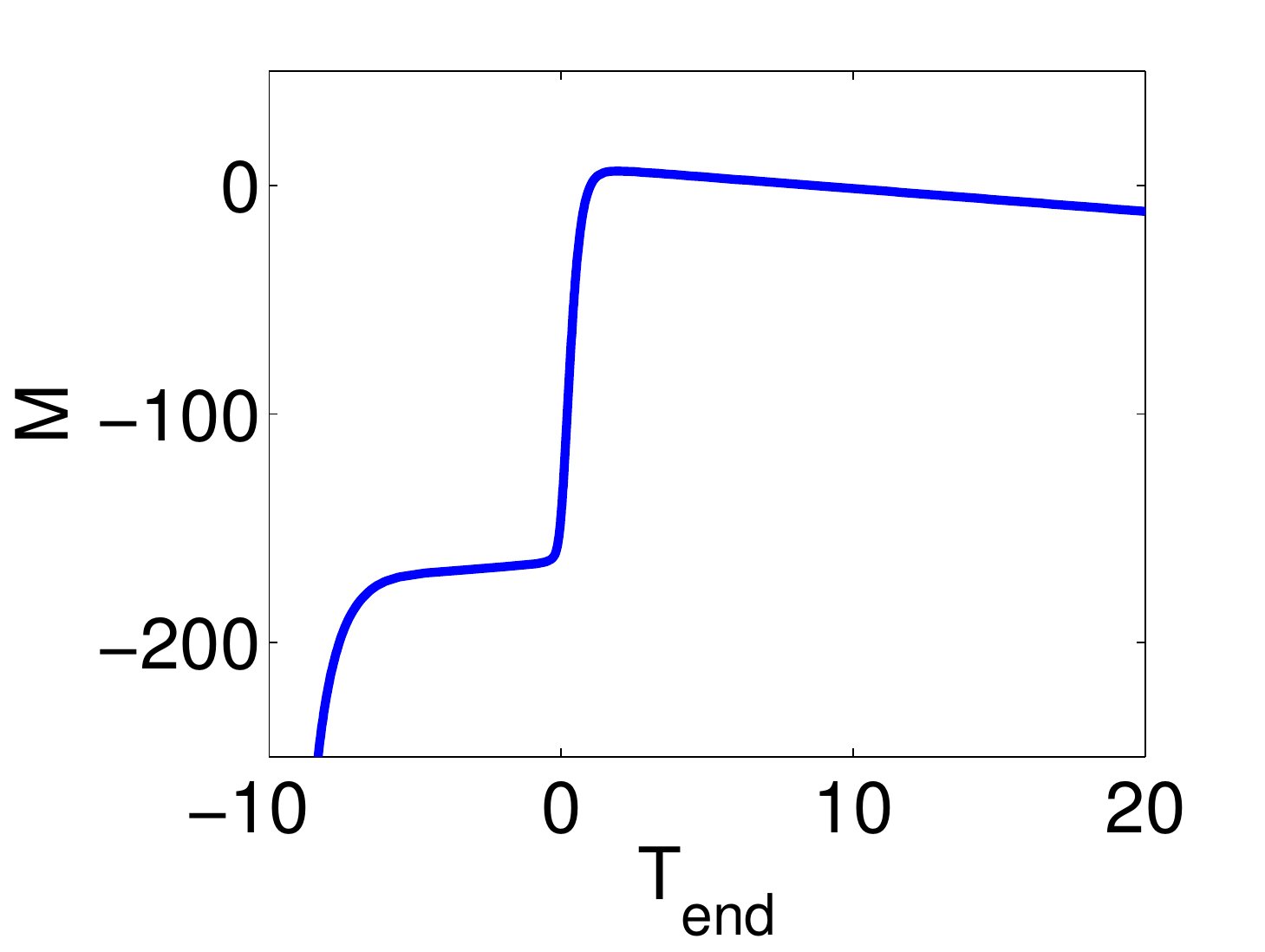}}        
        ~ 
        \caption{(a) Optimal path for $T_{\mathrm{end}} = 20$. (b) Plot of function $M$ that needs to be maximized w.r.t. $T_{\mathrm{end}}$. $t_{0} = -10$, $\epsilon = 1.25$, $D = 0.008$.}\label{Tend continuations}
\end{figure}

We would like to point out that the BVP \eqref{2nd
  order ODE} used to calculate the locally optimal path is valid for a
scalar time-dependent system and for finite (non-small) noise variance
$2D$. In the small noise limit, one can use minimum action methods to
find the optimal path, which can be applied to multiple dimensions
\cite{ren2004minimum}. Furthermore, according to
\citet{ren2004minimum}, in gradient systems, over an infinite time
interval, the optimal path becomes a minimum energy (where `energy' refers to the functional $F$ that is optimized) path that forms a
heteroclinic orbit between the two local minima of the
potential. However, Figure~\ref{Tend continuations}
demonstrates that even for relatively small noise levels $D$ (such as
$D=0.008$ as chosen for previous illustrations) we are far away from
the small noise limit such that $T_{\mathrm{end}}$ is of order $1$:
Figure~\ref{Tend continuation optimal path} shows the (locally)
optimal path $x(t)$ for $T_{\mathrm{end}} = 20$. We observe that for a
long time ($1<t<18$) the path $x(t)$ stays close to the saddle $U_+$
before eventually escaping to the chosen $x_{T} = 4$. As
Figure~\ref{Tend continuation 1} shows, the lingering of $x(t)$ close to
the saddle is only optimal for the fixed large $T_\mathrm{end}=20$. The functional $M = \log(F)$
increases for decreasing $T_\mathrm{end}$.

This implies that for positive (even small) noise variance $2D$ the
functional $F$ should also be optimized with respect to the traveling
time $T_\mathrm{end}$ of the path. We formulate the extended BVP
corresponding to critical points with respect to path and traveling
time in rescaled time on the base interval $[0,1]$. The BVP will then
be solved with standard continuation software AUTO \cite{D07}. The BVP \eqref{2nd order ODE}, rescaled
to $[0,1]$ is (split into two components):
\begin{align}
\label{scaled ODE1}
  \dot{x}_{1} &= x_{2}(T_{\mathrm{end}} - t_{0})\mbox{,}&
  x_{1}(0)&= x_{0}\mbox{,}\\
   \dot{x}_{2} &= g(x_{1},t)(T_{\mathrm{end}} - t_{0})\mbox{,} &
   x_{1}(1) &= x_{T}
\label{scaled ODE2}
\end{align}
where $t_{0}$ (fixed) and $T_{\mathrm{end}}$ (free) are the start and
end $t$ values and
\begin{eqnarray*}
g(x_{1},t) = 2D\dfrac{\partial V_{s}}{\partial x}(x_1(t),t)\mbox{.}
\end{eqnarray*}
The solution of \eqref{scaled ODE1}--\eqref{scaled ODE2} is a critical
point of $F$, given in \eqref{Functional}, among all possible paths
connecting from $x_{0}$ to $x_{T}$ in a fixed time $T =
T_{\mathrm{end}} - t_{0}$.  The function $M = \log(F)$, written for
the rescaled path is
\begin{align}
\label{Integral M}
&M = \int\limits_{0}^{1}\dfrac{U_{0} - U_{T}}{2D} - \bigg[\dfrac{x_{2}(t)^{2}}{4D}
+ V_{s}(x_{1}(t),t)\bigg](T_{\mathrm{end}} - t_{0})\mathrm{d}t
\\
\nonumber
&\mbox{where}\quad U_{0} = U(x_{0},t_{0})\mbox{,}\quad
U_{T} = U(x_{T},T_{\mathrm{end}})\mbox{.}
\end{align}
Paths maximizing $F$ (and, hence, $M$) also maximize the probability
of realizations of SDE \eqref{SDE} following it.  Figure \ref{Tend
  continuation 1} plots $M$ along paths satisfying \eqref{scaled ODE1}--\eqref{scaled ODE2}
for a range of end times $T_{\mathrm{end}}$. Its maximum corresponds
to the time $T_\mathrm{end}$ for which the functional $M$ is (locally)
maximal among the range of $T_\mathrm{end}$ shown (this is for the
fixed positive but small noise variance $2D=0.016$). We now extend the
BVP \eqref{scaled ODE1}--\eqref{scaled ODE2} to include the criticality of $T_\mathrm{end}$
into the optimization problem. We outline the BVP of the variational
problem for the general case here. The specific case for this example
is in Appendix~\ref{sec:bvpderiv}.

Introducing the derivatives of $x_{1}$ and $x_{2}$ w.r.t. $T_{\mathrm{end}}$ as
\begin{align*}
z_{1}(t) = \dfrac{\partial x_{1}(t)}{\partial T_{\mathrm{end}}}, \qquad z_{2}(t) = \dfrac{\partial x_{2}(t)}{\partial T_{\mathrm{end}}}\mbox{,}
\end{align*}
these derivatives satisfy
\begin{eqnarray}
\label{z ODE}
\begin{aligned}
  \dot{z}_{1} &= x_{2} + z_{2}(T_{\mathrm{end}} - t_{0})\mbox{,}&
  z_{1}(0)&= 0\mbox{,} \\
  \dot{z}_{2} &= g(x_{1},t) + \dfrac{\partial g(x_{1},t)}{\partial
    x_{1}}z_{1}(T_{\mathrm{end}} - t_{0})\mbox{,} &
   z_{1}(1) &= 0\mbox{.}
\end{aligned}
\end{eqnarray}
Critical points of $M(x_{0},t_{0},x_{T},T_{\mathrm{end}},x_{1}(\cdot),x_{2}(\cdot))$, given in equation (\ref{Integral M}), w.r.t. $T_{\mathrm{end}}$ satisfy
\begin{eqnarray*}
\dfrac{\partial M}{\partial T_{\mathrm{end}}} +  \dfrac{\partial M}{\partial x_{1}}\dfrac{\partial x_{1}}{\partial T_{\mathrm{end}}} + \dfrac{\partial M}{\partial x_{2}}\dfrac{\partial x_{2}}{\partial T_{\mathrm{end}}} = 0\mbox{,}
\end{eqnarray*} 
which produces the integral condition:
\begin{multline}
0=m:=\int\limits_{0}^{1}\bigg[\dfrac{1}{2D}\dfrac{\partial U(x_{T},T_{\mathrm{end}})}{\partial T_{\mathrm{end}}} + \dfrac{x_{2}(t)^{2}}{4D} + V_{s}(x_{1}(t),t) \\
+ \bigg(\dfrac{x_{2}(t)z_{2}(t) + g(x_{1}(t),t)z_{1}(t)}{2D}\bigg)(T_{\mathrm{end}} - t_{0})\bigg]\mathrm{d}t
\label{Integral condition}
\end{multline} 
Therefore, have to solve the four-dimensional BVP \eqref{scaled ODE1},
\eqref{scaled ODE2}, \eqref{z ODE} for $x_{1}(t)$, $x_{2}(t)$,
$z_{1}(t)$, $z_{2}(t)$ with the additional integral condition
\eqref{Integral condition} and the additional parameter
$T_{\mathrm{end}}$. We use AUTO (Version: AUTO-07P) \cite{AUTOweb} to study the solutions of
\eqref{scaled ODE1}, \eqref{scaled ODE2}, \eqref{z ODE},
\eqref{Integral condition} in dependence of the system parameters $D$
and $\epsilon$.

\section{Sequence of continuation steps for the optimal path to escape in optimal time}
\label{sec:steps}
Since \eqref{scaled ODE1}, \eqref{scaled ODE2}, \eqref{z ODE},
\eqref{Integral condition} is nonlinear we need a sequence of
initialization steps to arrive at the optimal path for particular
desired values of ramping speed parameter (initially $\epsilon =
1.25$, close to critical value $\epsilon_c=4/3$) and noise variance
(initially $2D = 0.1$). An advantage of using continuation is that
once we have obtained an optimal path in an optimal time for a
particular set of parameters we are free to perform a systematic
parameter study of solutions of \eqref{scaled ODE1}, \eqref{scaled
  ODE2}, \eqref{z ODE}, \eqref{Integral condition} varying the
ramping speed $\epsilon$ and noise level $D$.

\subsection{List of free parameters} First we discuss some of the parameters used and reasoning for their initial values, as given in Table \ref{Table of Parameters}.
\begin{table}[ht]
\caption{\label{Table of Parameters}Types of parameters used in continuation steps and their initial values}
\begin{ruledtabular}
\begin{tabular}{p{2.05cm}p{2cm}p{2cm}p{2cm}}
    System (fixed) parameters & Continuation parameters & Bifurcation parameters & Monitoring parameter\\[0.5ex] \hline
    {$\!\begin{aligned}\ \\[-3ex]  
               p &= 1 \\    
               \lambda_{\max} &= 3 \\
               t_{0} &= -10 \\
               x_{0} &= -1 \end{aligned}$} &  {$\!\begin{aligned}\ \\[-3ex]    
               T_{\mathrm{init}} &= 0 \\    
               x_{T} &= x_{0} \\
               T_{\mathrm{end}} &= -9 \\
               m &= m_{0} \end{aligned}$} &  {$\!\begin{aligned} \ \\[-3ex]   
               \epsilon &= 1.25 \\    
               D &= 0.05 \\ \\ \\
                \end{aligned}$} & {$\!\begin{aligned}\ \\[-3ex]    
               M &= M_{0} \\ \\ \\ \\  \end{aligned}$} \\
    \end{tabular}
    \end{ruledtabular}
\end{table}
We introduce the factor $T_{\mathrm{init}}$ in equation \eqref{scaled
  ODE2} as an artificial parameter. Thus, \eqref{scaled
  ODE1}--\eqref{scaled ODE2} changes to
\begin{eqnarray*}
\dot{x}_{1} &&= x_{2}(T_{\mathrm{end}} - t_{0})\\
\dot{x}_{2} &&= g(x_{1},t)(T_{\mathrm{end}} - t_{0})T_{\mathrm{init}}\mbox{,}
\end{eqnarray*}
giving a trivial system ($\dot x_1=x_2(T_{\mathrm{end}} - t_{0})$,
$\dot{x}_{2}=0$) for $T_\mathrm{init}=0$, connecting it to system
\eqref{scaled ODE1}--\eqref{scaled ODE2} via a continuation in
$T_\mathrm{init}$ from $0$ to $1$. Furthermore, we initially choose
$x_{T} = x_{0}$ and $T_{\mathrm{end}} = t_0+1=-9$, close to the
initial time value $t_{0}$. Finally, the parameters $M$ and $m$ are
used to monitor the values of the integrals in equations
\eqref{Integral M} and \eqref{Integral condition}, respectively, such
that sign changes of $m$ correspond to critical values of the
functional $M$ (and, hence, $F$). The initial values of $M$ and $m$
are set equal to the integrals in \eqref{Integral M} and
\eqref{Integral condition} along the initial
path. 
Outlined in Table \ref{Continuation steps} is a brief summary of the
continuation steps performed to create an optimal path in an optimal
time. We proceed with a brief explanation for each of the continuation
steps with a more in depth discussion provided in Appendix \ref{app:continuation steps}.
\begin{table}[ht]
\caption{\label{Continuation steps}Summary of continuation steps to perform in order to achieve an optimal path for escape in an optimal time (in brackets are the values used)}
\begin{ruledtabular}
    \begin{tabular}{ p{0.65cm}  p{1.8cm}  l  l  p{1.5cm} }
    Step $\#$ & Continuation parameter & Initial value & End value & Other free parameters \\ \hline
    1 & $T_{\mathrm{init}}$ & $(0)$ & $(1)$ & $m$, $M$ \\ 
    2 & $x_{T}$ & $x_{0}$ $(-1)$ & $\gg 1$ $(4)$ & $m$, $M$ \\
    3 & $T_{\mathrm{end}}$ & $\sim t_{0}$ $(-9)$ & $\gg 1$ $(20)$ & $m$, $M$ \\ 
    \end{tabular}
    \end{ruledtabular}
\end{table}
\subsection{Step 1: continuation of $T_{\mathrm{init}}$ from $0$
  to $1$} The first
step 
ends in a solution of the full system of equations \eqref{scaled
  ODE1}-\eqref{scaled ODE2}. This is the short orbit (blue) from
$x_1=S_-=-1$ to $x_1=-1$ in Figure~\ref{PP sketch}. The parameters $M$
and $m$ are kept free during the continuation, monitoring the integrals
in \eqref{Integral M} and \eqref{Integral condition}.

\subsection{Step 2: continuation in $x_{T}$} 
\begin{figure}[ht]
        \centering       
                \includegraphics[scale = 0.4]{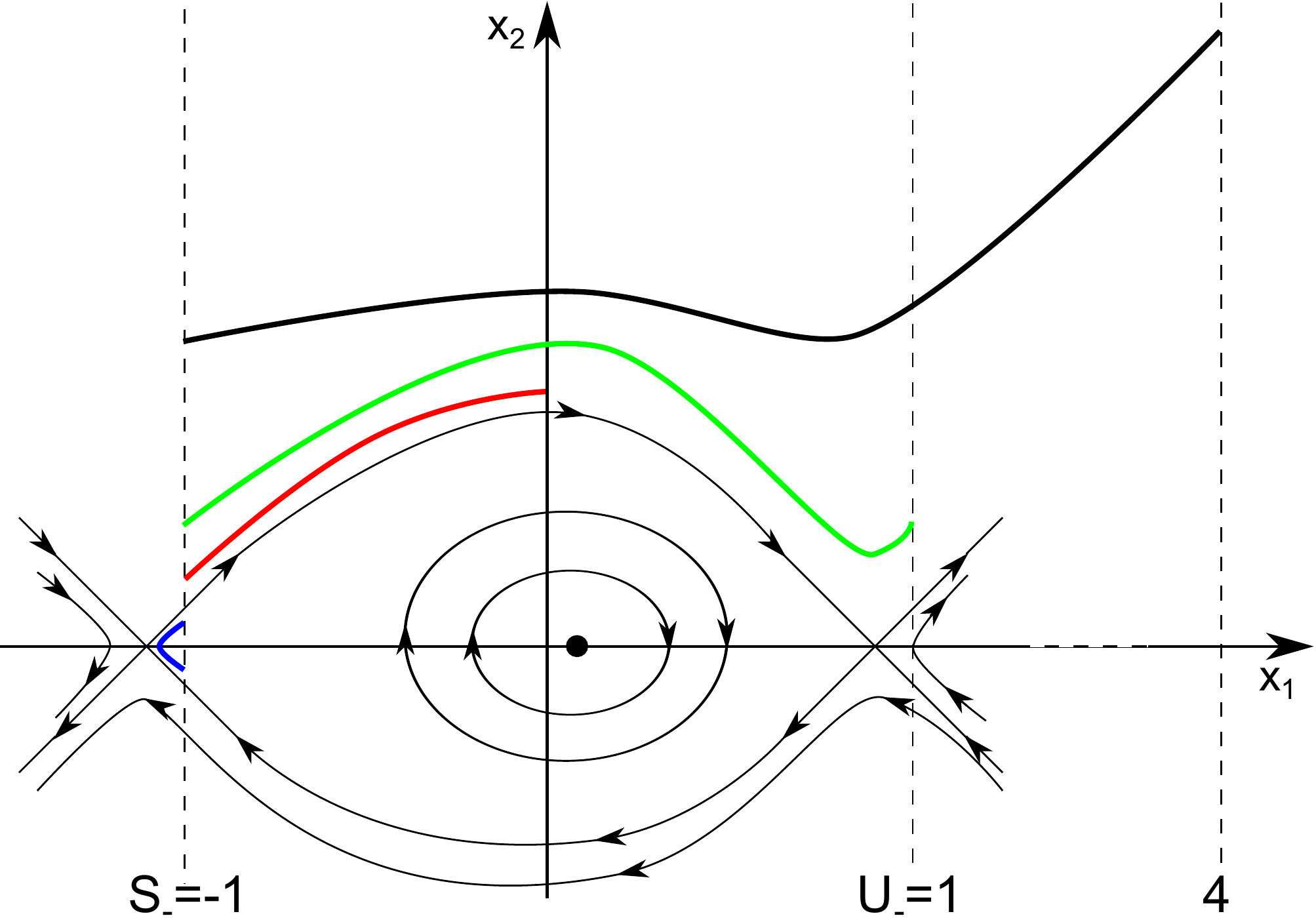}
                \caption{Illustration of trajectories (colored) at different stages of the $x_{T}$ continuation step superimposed on the phase portrait for the full rate-induced system, when stationary, $\lambda = 0$. $x_{T} = x_{0} = -1$ (blue), $x_{T} = 0$ (red), $x_{T} = 1$ (green), $x_{T} = 4$ (black).}
                \label{PP sketch}
\end{figure}

The result of Step~1 is a
path maximizing $M$ when traveling from $x_{0}=-1$ to $x_T=-1$ in
unit time ($t_{0} = -10$, $T_{\mathrm{end}} = -9$), where $\lambda(t)$
is close to stationary ($\lambda\approx0$). Step~2 changes the right
boundary value $x_{T}$ to the desired location. Guided by our aim to
find most likely paths for escape, we perform a continuation to $x_{T}
= 4$ (we count a trajectory of \eqref{SDE} that reaches $x_T=4$ as
having escaped). Figure~\ref{PP sketch} shows a sequence of solutions
of \eqref{scaled ODE1}--\eqref{scaled ODE2} (colored) for the
continuation stages of $x_{T}$ superimposed onto the phase portrait of the system
\eqref{scaled ODE1}--\eqref{scaled ODE2} for $\lambda = 0$.
The results of this step is a path maximizing $M$ that connects
$x_0=-1$ and $x_T=4$ in a very short time period
($T_\mathrm{end}-t_0=1$, see Figure \ref{Optimal path after xT
  continuation} in Appendix \ref{app:xT continuation}). Clearly this is not the optimal time to
make this transition and so the next step is to continue in
$T_{\mathrm{end}}$ to get to more realistic timings of escape.
   
\subsection{Step 3: continuation in $T_{\mathrm{end}}$} We increase
$T_{\mathrm{end}}$ to a large value, monitoring $M$ and $m$. Figure~\ref{Tend
  continuation 1} shows the graph of $M$ over $T_\mathrm{end}$, which
has a pronounced maximum at
$T_\mathrm{end}\approx1.5$. Critical points of $M$ are detected when $m$ changes sign, and therefore gives the optimal path for an optimal time provided $M$ is a maximum. 


The optimal path constructed by the above steps will be
systematically continued in the system parameters $D$ (noise variance $2D$) and
$\epsilon$ (ramping speed of $\lambda(t)$) in Section~\ref{sec:timing}. 

\section{Optimal path for escape for noise and rate-induced tipping}
\label{sec:optimal}
\begin{figure}[ht]
        \centering
        \subcaptionbox{Time profile.\label{TP optimal escape}}[0.45\linewidth]
                {\includegraphics[scale = 0.3]{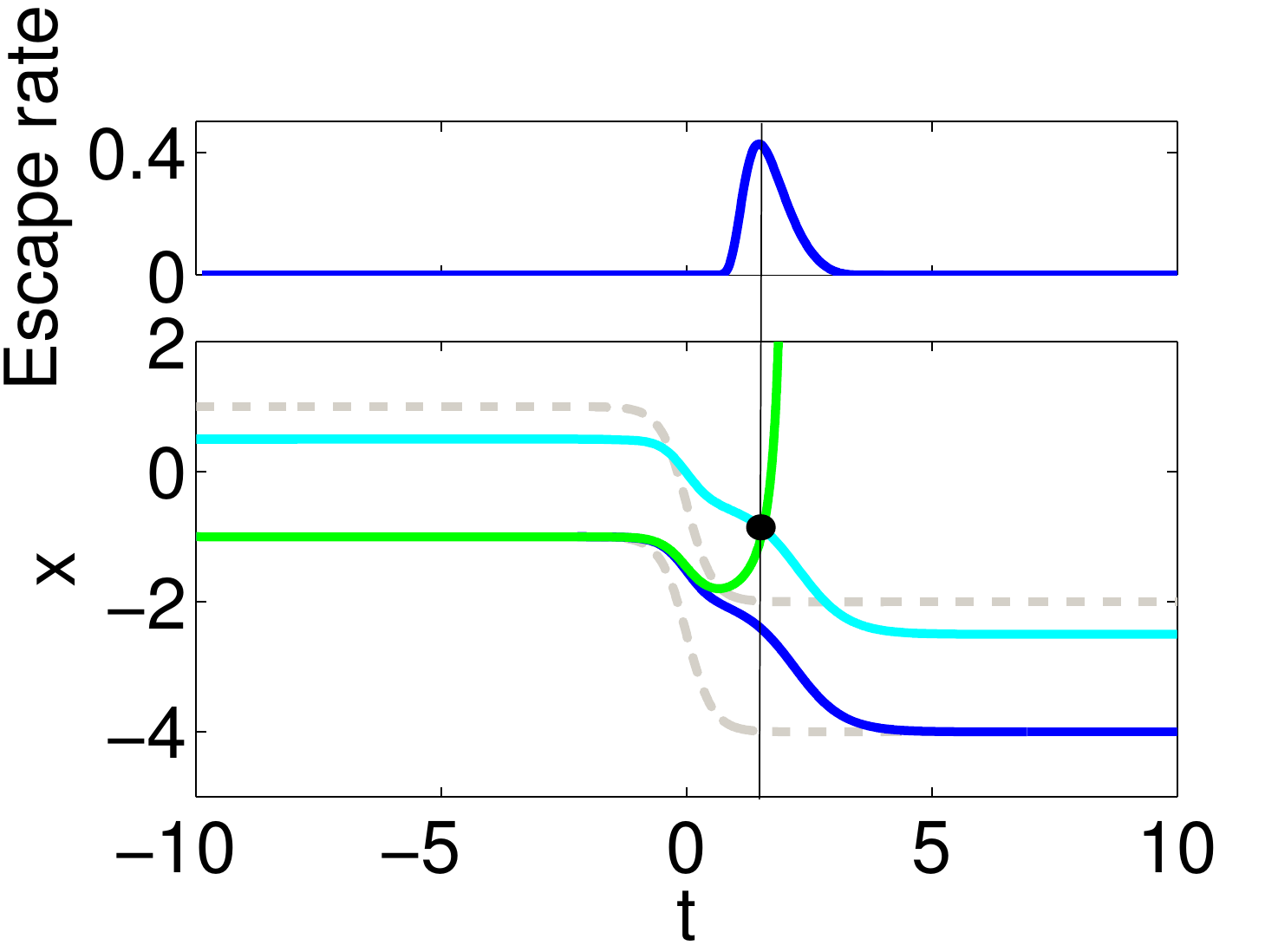}}
        \hfill 
        \subcaptionbox{Phase plane.\label{PP optimal escape}}[0.45\linewidth]
                {\includegraphics[scale = 0.3]{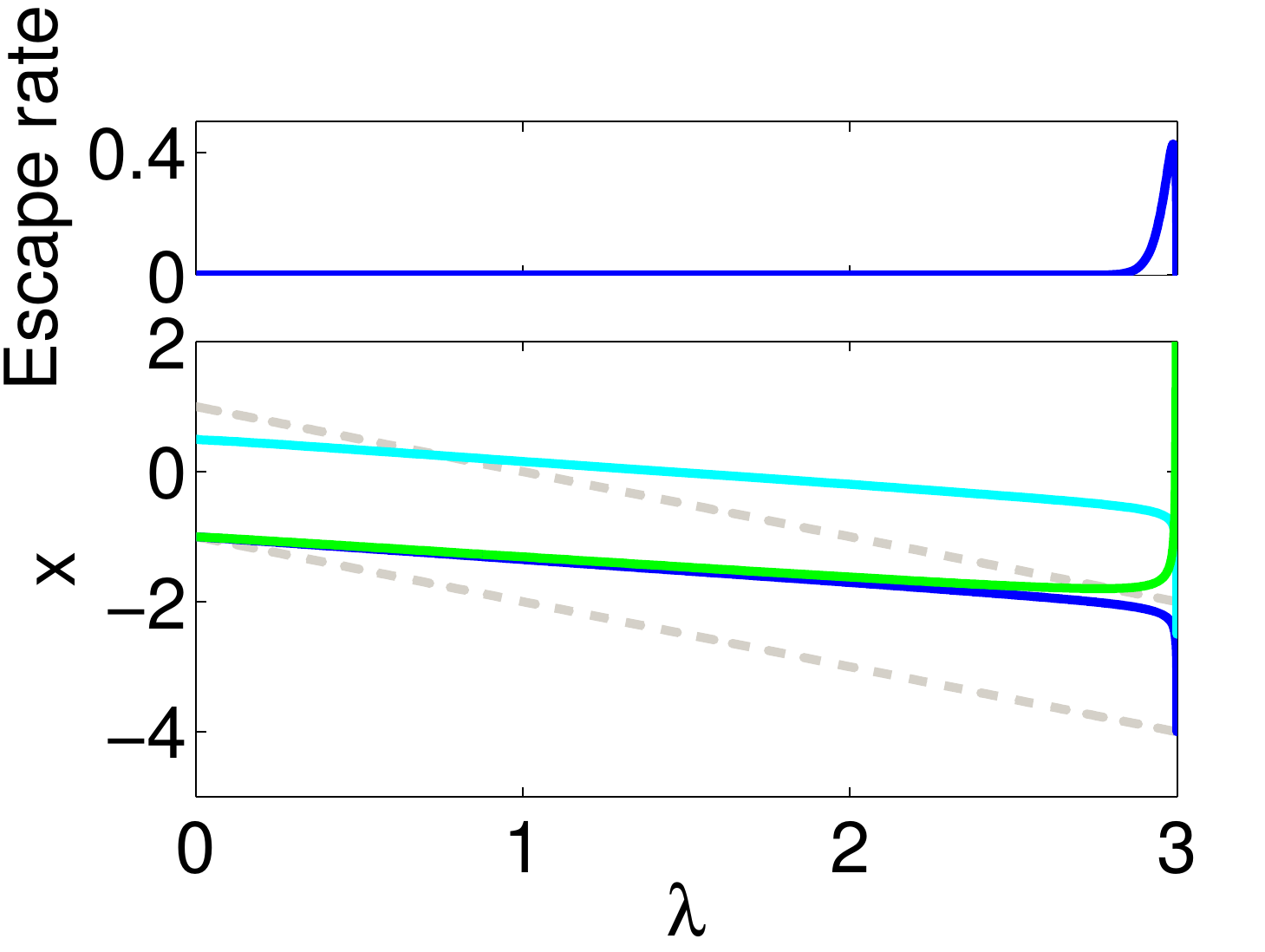}}
        ~ 
        \caption{Time profile and phase plane with escape rate from simulations (top) and optimal path added in green for $\epsilon = 1.25$ and $D = 0.008$. Dashed curves represent stable $W^{s}_{0}$ and unstable $W^{u}_{0}$ quasi-equilibria, the dark blue curve is the unstable manifold $W^{u}(S_{-})$, the bright blue curve is the threshold curve $\tilde{x}(t)=x^u(t)+y$ with $y = 1.5$.}\label{Escape rate optimal}
\end{figure}

This section compares the optimal path for escape with the escape calculated directly from the solutions of
    \eqref{FPE}. We are interested in the timing of escape, defined as the timing of
crossing certain threshold curves. We do not want to have to rely on
running full Monte Carlo simulations, but instead use the optimal path
theory developed by \citet{chaichian2001path,zhang2008theory}. We include the optimal path for escape (in
green) into Figure~\ref{Escape rate} to compare the most likely timing
of escape (the peak in the escape rate, measured at the threshold $\tilde x$) with
the time $t$ the threshold $\tilde{x}$ intersects with the optimal path as computed through continuation, see Figure
\ref{Escape rate optimal}.
The corresponding phase portrait  in Figure \ref{PP optimal escape} shows that
both the simulations and optimal path suggest that the escape does not
happen until just short of $\lambda = 3$, the moment the potential
well or steady states $W^{s}_{0}$, $W^{u}_{0}$ are coming to a rest.

Figure
\ref{Escape rate optimal} illustrates that the optimal path matches
the mode (peak) of the escape rate well. In general, if the escape
rate over time is unimodal with a sharp peak then the time profile of
the optimal path is a good description of this peak. More precisely, the mode of the escape rate
occurs very close to a time $t$ for which $\tilde x= x_1(t)$, where
$x_1$ is the first component of the optimal path, the solution of the
extended BVP \eqref{scaled ODE1}, \eqref{scaled ODE2}, \eqref{z ODE},
\eqref{Integral condition}.

\subsection{Dependence on choice of threshold curve $\tilde{x}(t)$:}
The choice of the threshold curve $\tilde{x}(t)$ in Figure~\ref{Escape
  rate optimal} is at $\tilde x(t)=x^u(t)+y$ with $y=1.5$ (recall that $x^u(t)$ is the trajectory of the deterministic part of  \eqref{SDE}. Figure \ref{Escape rate color
  optimal} shows a color plot for the escape rate at the
threshold depending on the distance $y$. The distance of the optimal
path from the unique trajectory $x^{u}(t)$ added in white. This
highlights that, provided the threshold $\tilde{x}(t)$ is sufficiently
far from $x^{u}(t)$, the optimal path will cross the threshold at the
same moment the escape rate through the threshold is at its peak.

\begin{figure}[ht]
        \centering
                \includegraphics[scale = 0.3]{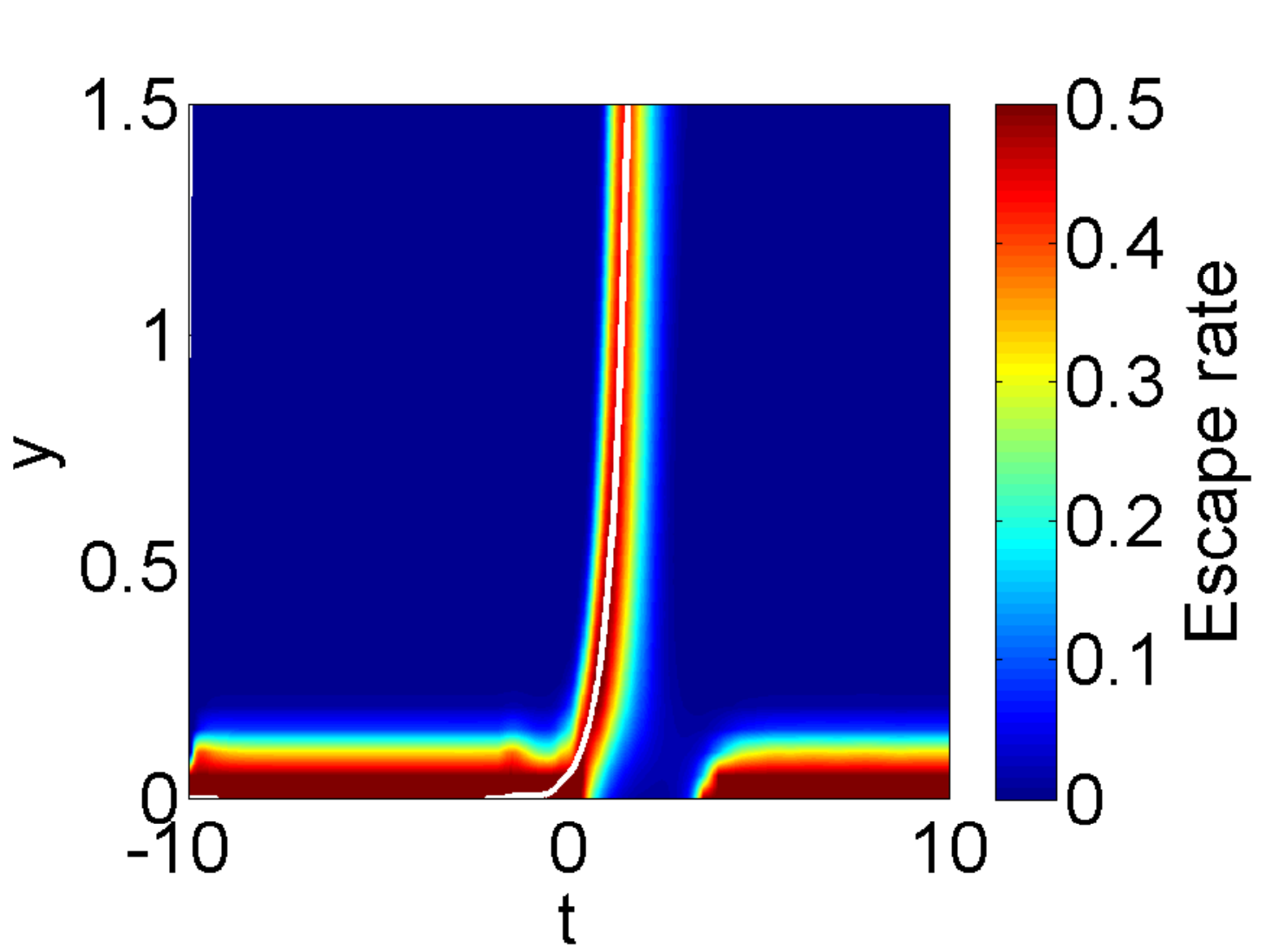}
                \caption{Color plot of the escape rate at the threshold $\tilde{x}(t)$ depending on the distance $y$ from the unstable manifold $W^{u}(S_{-})$. Distance of optimal path from $x^{u}(t)$ added in white, $\epsilon = 1.25$, $D = 0.008$.}
                \label{Escape rate color optimal}
\end{figure} 

\subsection{\change{Dependence on starting time $t_{0}$ and starting position $x_{0}$:}}
\change{We emphasize that the optimal path calculates the local
  optimum, i.e. assuming the system has not tipped before the ramping
  shift begins. We set $t_{0} = -10$, $x_{0} = -1$ to represent
  starting at the bottom of the potential well at time
  $-\infty$. Changing the starting time, for example to $t_{0} = -15$
  or $t_{0} = -5$ has no effect on when the optimal path escapes and
  only extends or shortens the time profile of the path presented in
  Figure \ref{Escape rate optimal}(a). Likewise, changing the starting
  position, provided it is still inside the well, has no effect on
  when the optimal path escapes. For slightly different $x_0$ the
  optimal path will converge onto the path in Figure \ref{Escape rate
    optimal}(a) before $t=0$ and then follow the same path for
  escape.}

\begin{figure}[ht]
        \centering
        \subcaptionbox{Time profile.\label{TP density}}[0.45\linewidth]
                {\includegraphics[scale = 0.3]{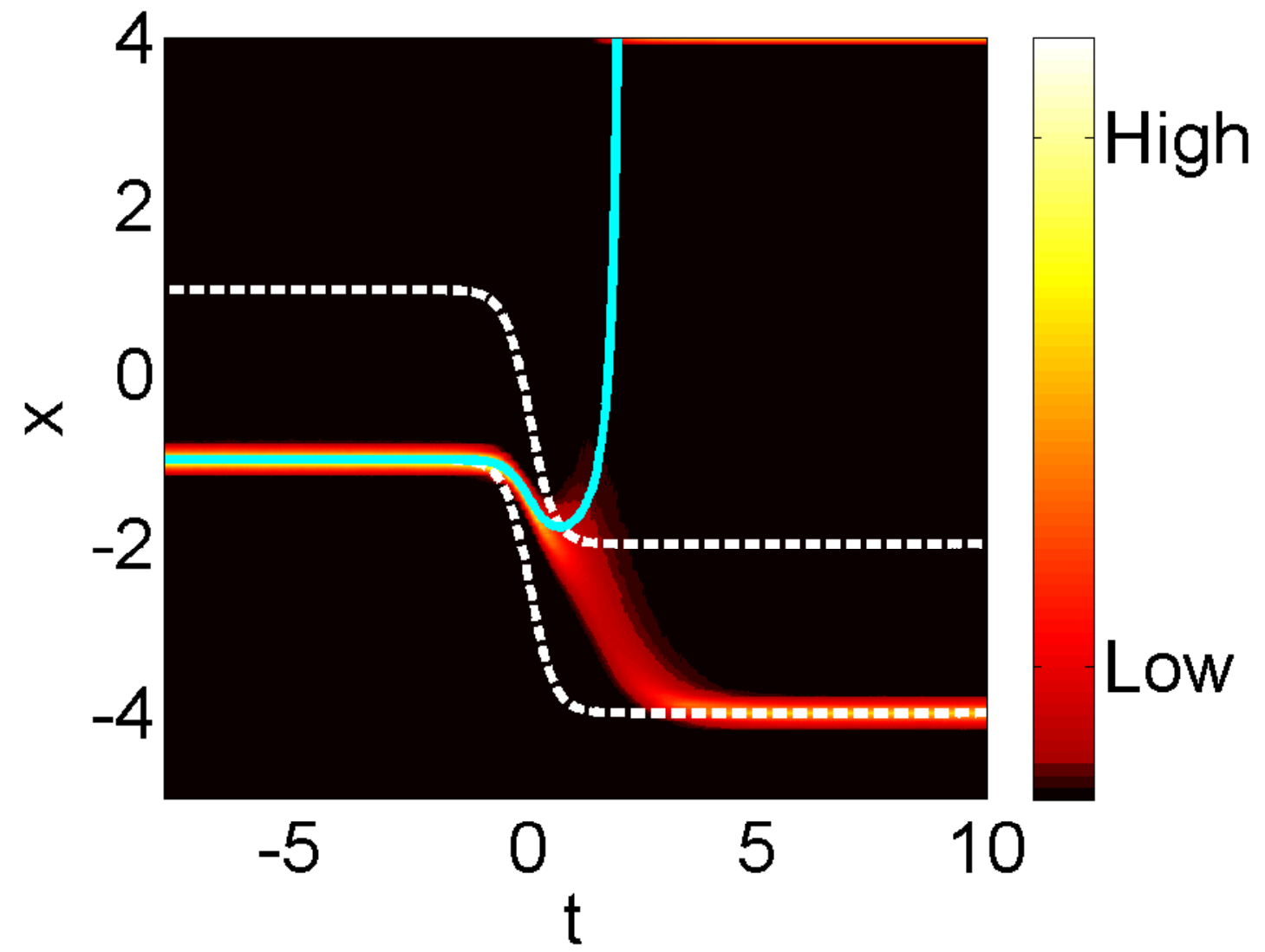}}
        \hfill 
        \subcaptionbox{Phase plane.\label{PP density}}[0.45\linewidth]
                {\includegraphics[scale = 0.3]{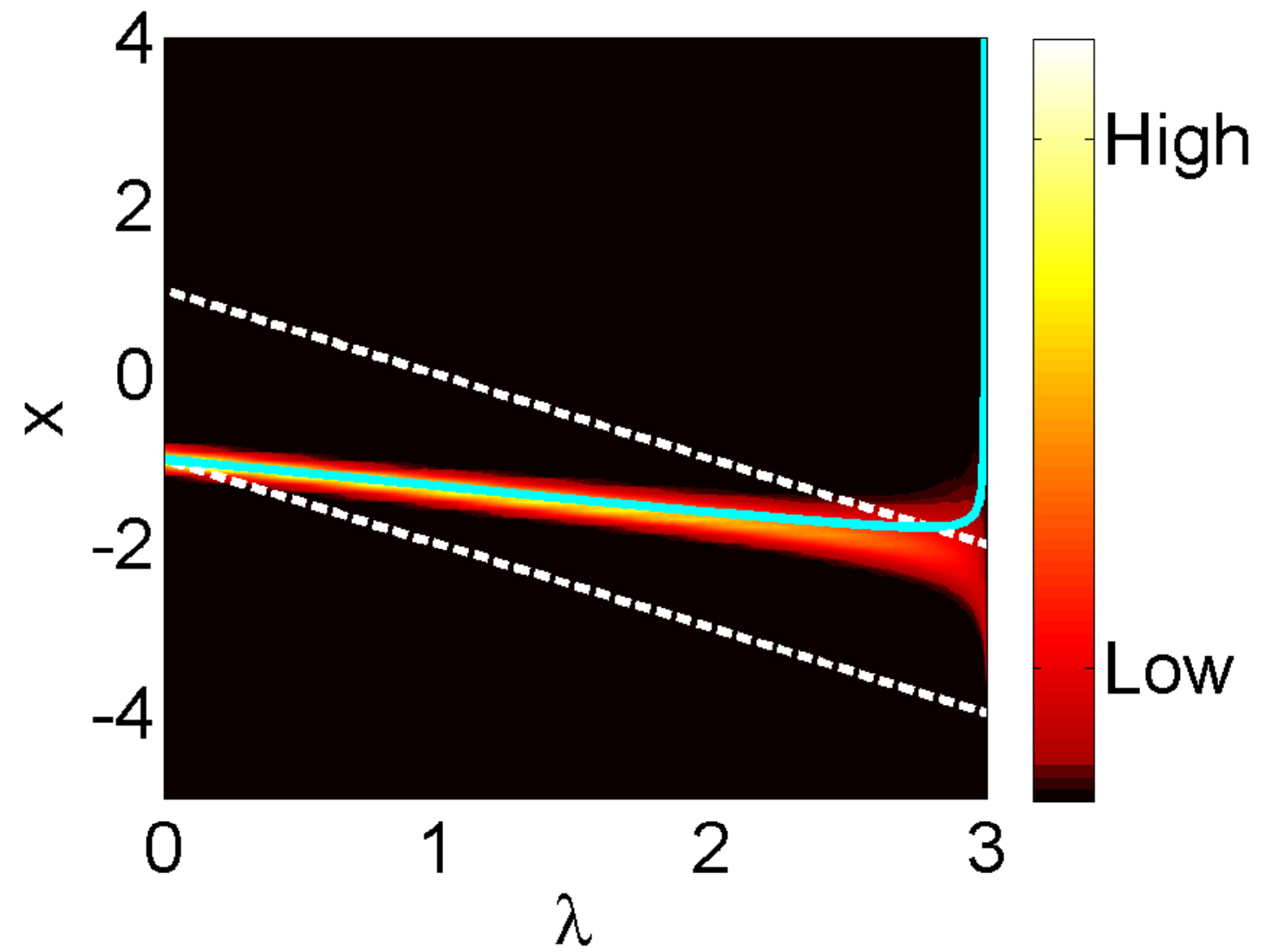}}
        ~ 
        \caption{Time profile and phase plane for density plot of simulations and optimal path added in bright blue for $\epsilon = 1.25$ and $D = 0.008$. Dashed white curves represent stable and unstable equilibria.}\label{Density}
\end{figure}

Figure~\ref{Density} shows how the $x_1(t)$ component of the optimal path is connected to the evolution of the density of realizations.

In the time profile plot, initially the spread of the distribution is
very narrow centered around the steady state $W^{s}_{0}$, due to the
small noise level $D$. When the system shifts, this distribution
widens reflected by a lower density over a larger $x$ range. Once the
shift stops, the density gradually becomes concentrated again, but
some realizations have escaped, indicated by the elevated density at
$x = 4$. Initially, the optimal path is right at the mean of the
distribution. The time when the optimal path deviates from the mean
equals the time the density in the simulation is at its widest and
where the additional mode (at $x=4$) appears. This once again suggests
that the optimal path, derived from BVP \eqref{scaled ODE1},
\eqref{scaled ODE2}, \eqref{z ODE}, \eqref{Integral condition},
describes the escape of realizations of the stochastic differential
equation \eqref{SDE}.

\section{Timing of escape in 2 parameter plane}
\label{sec:timing} 
\begin{figure}[ht]
        \centering
        \subcaptionbox{\label{Tend intersection}}[0.45\linewidth]
                {\includegraphics[scale = 0.3]{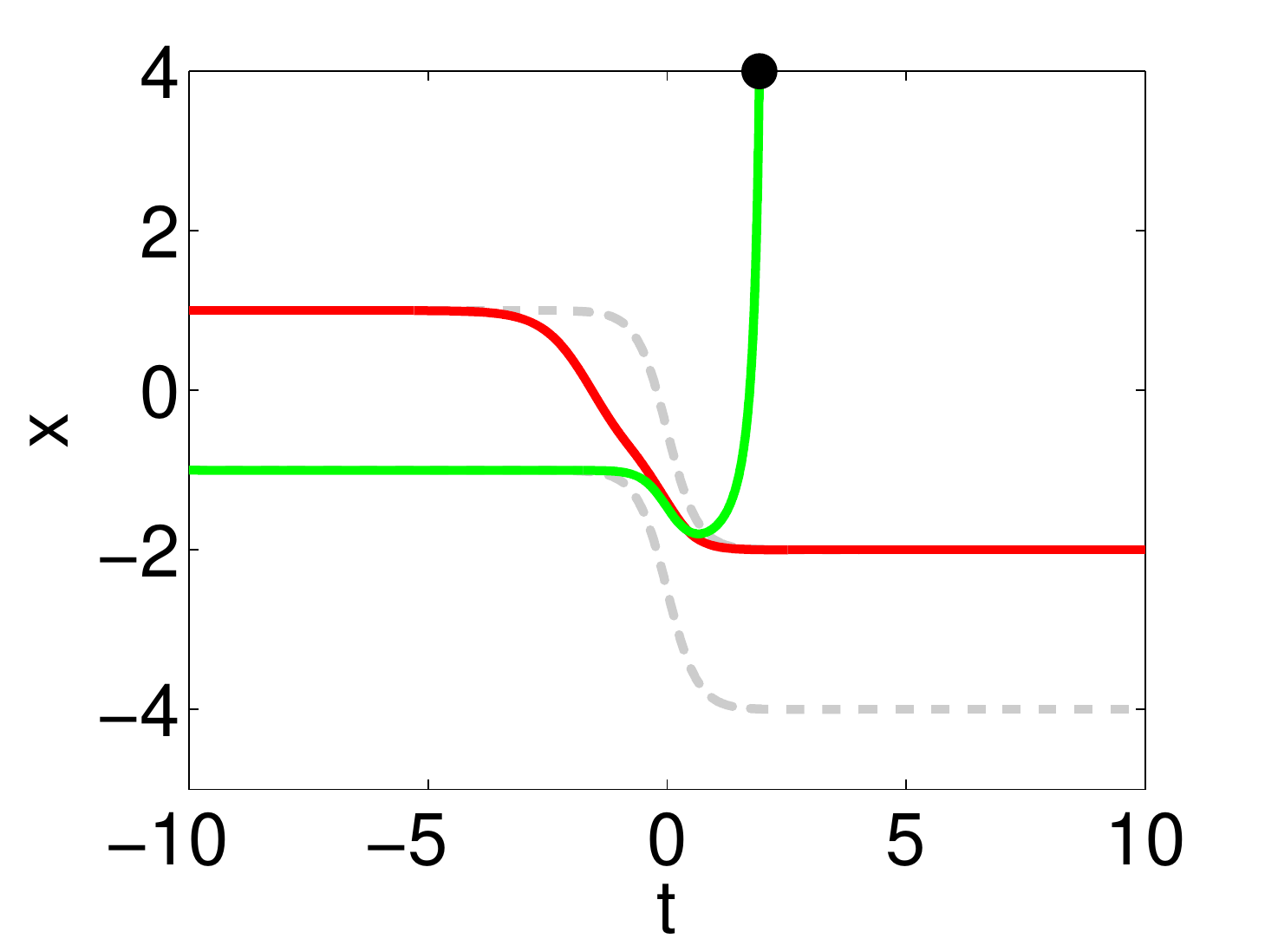}}
       \hfill 
        \subcaptionbox{\label{tcross intersection}}[0.45\linewidth]
                {\includegraphics[scale = 0.3]{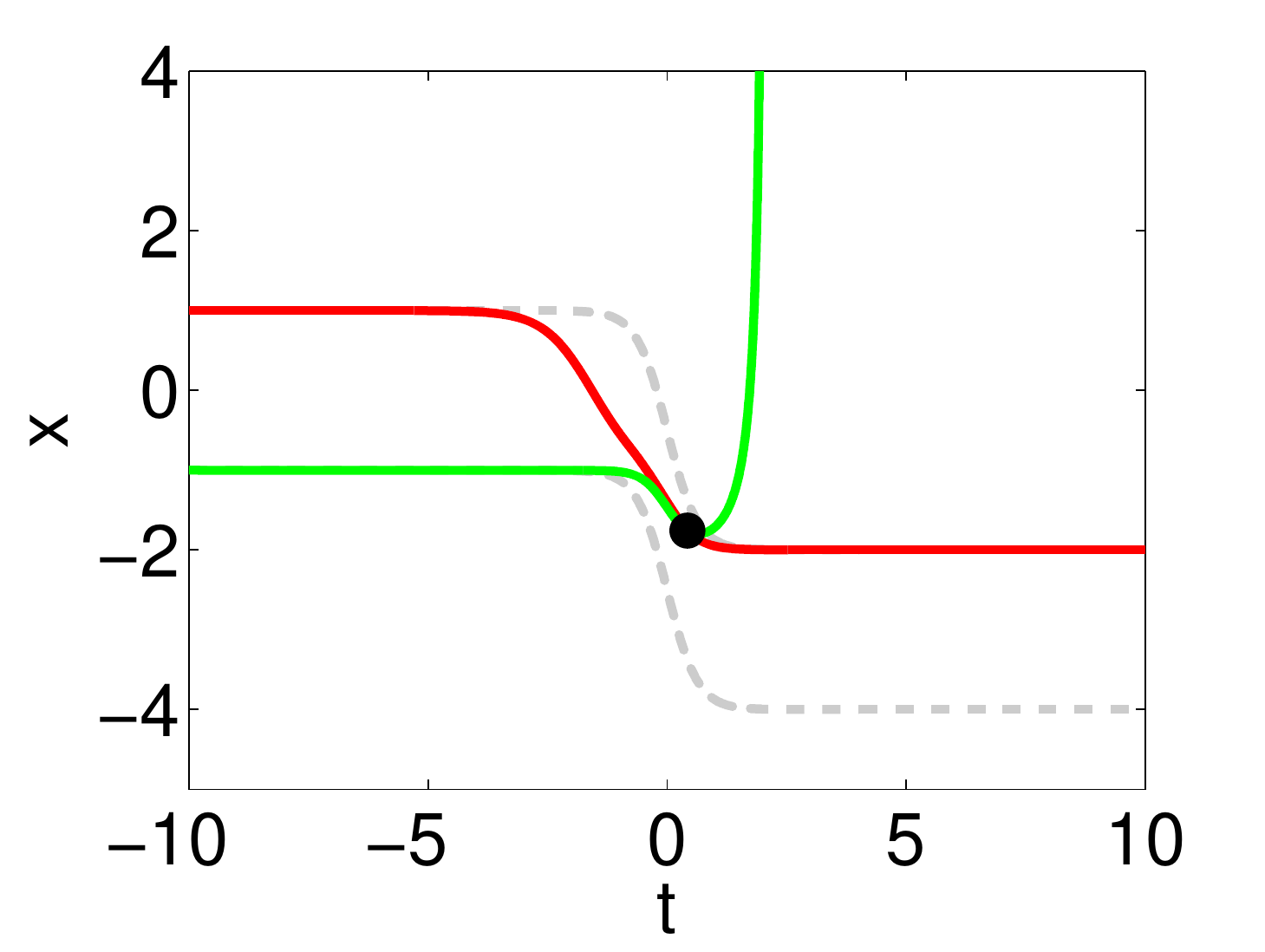}}
        ~ 
          \\
          \subcaptionbox{\label{Tend color}}[0.45\linewidth]
                {\includegraphics[scale = 0.3]{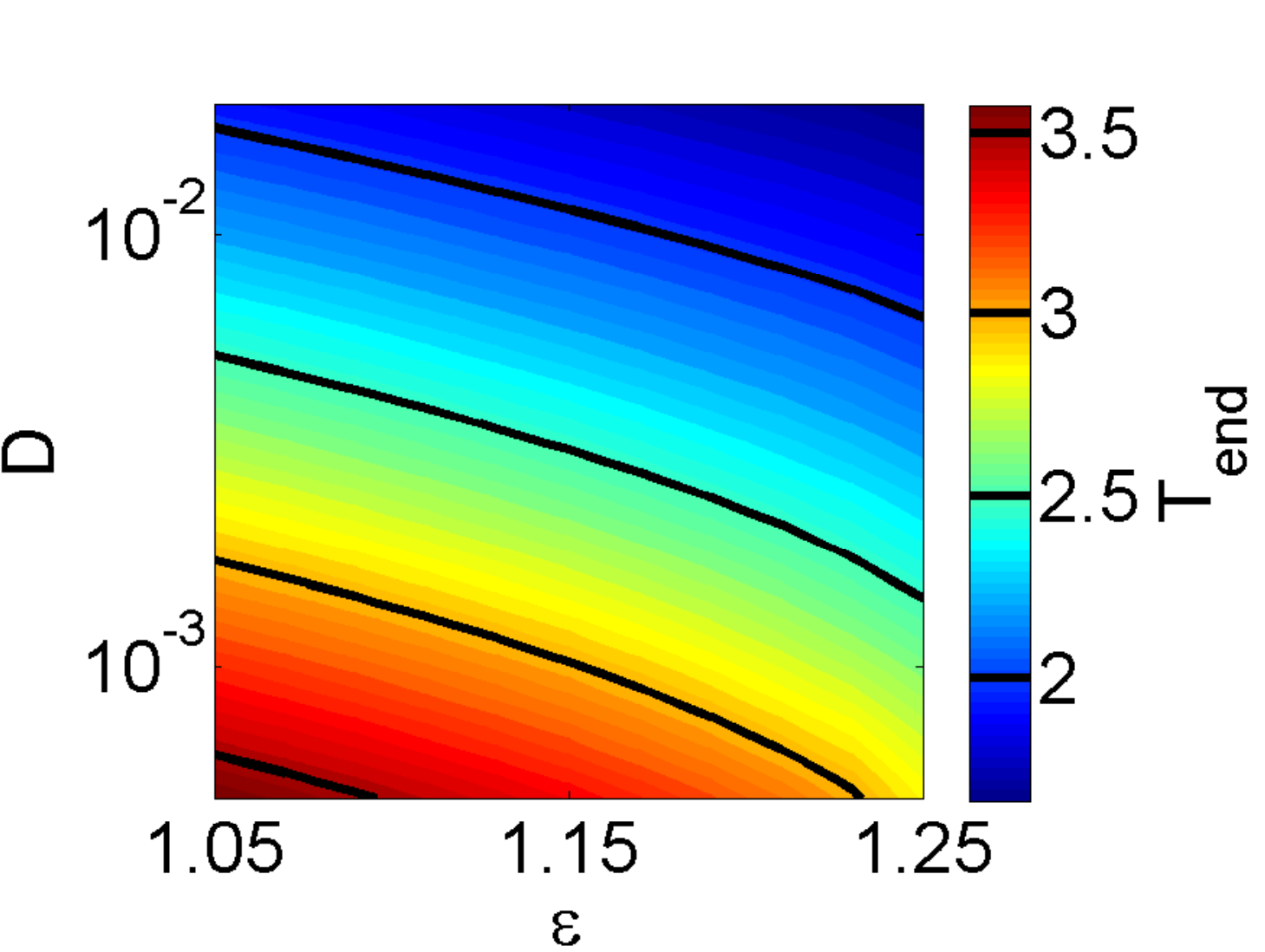}}
        \hfill 
        \subcaptionbox{\label{tcross color}}[0.45\linewidth]
                {\includegraphics[scale = 0.3]{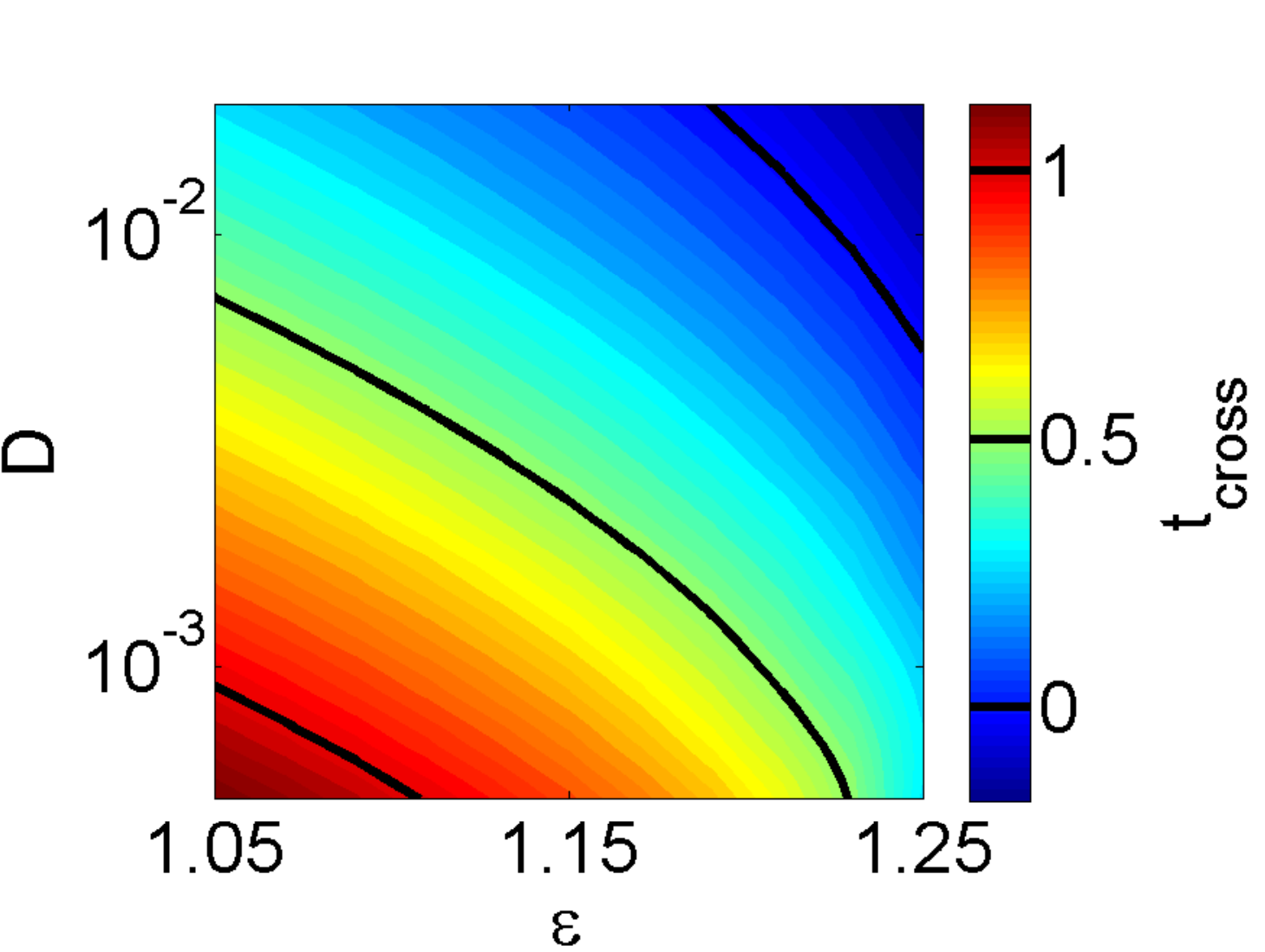}}
        ~ 
          \\
          \subcaptionbox{\label{Tend cross section}}[0.45\linewidth]
                {\includegraphics[scale = 0.3]{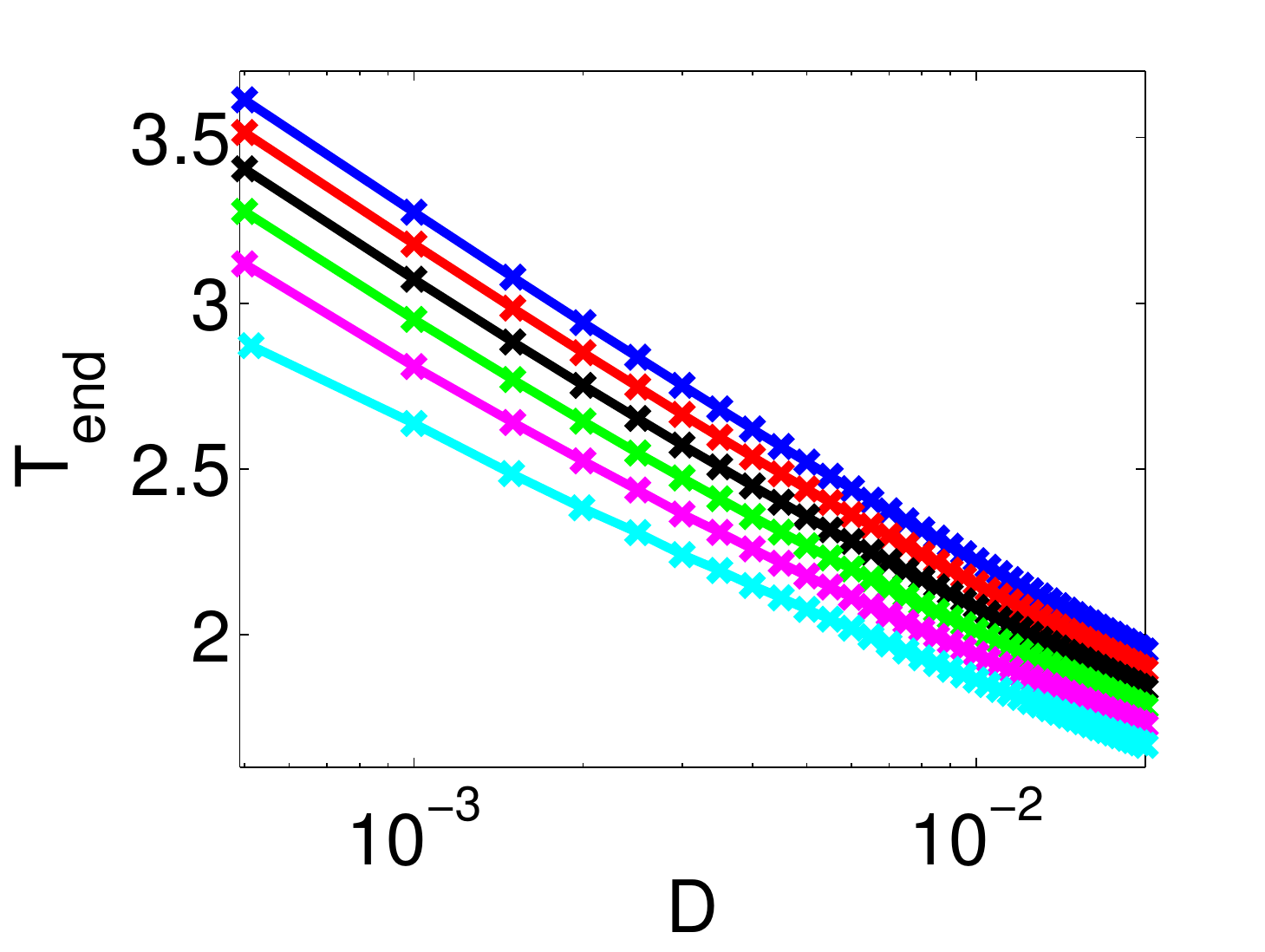}}
       \hfill 
        \subcaptionbox{\label{tcross cross section}}[0.45\linewidth]
                {\includegraphics[scale = 0.3]{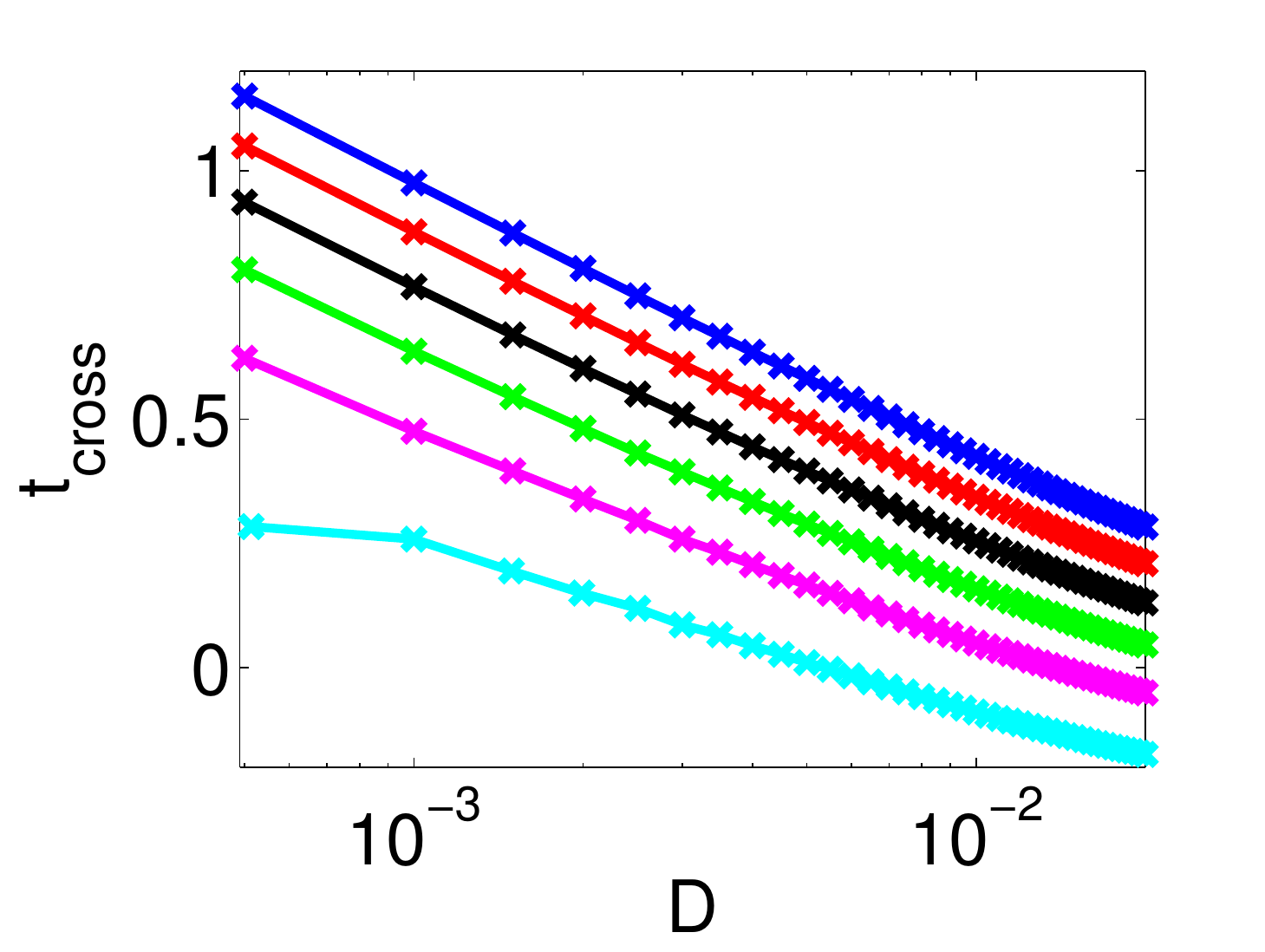}}
        ~ 
        \caption{Plots for end time $T_{\mathrm{end}}$ (a),(c),(e) and the crossing time $t_{\mathrm{cross}}$ with the stable manifold $W^{s}(U_{+})$ (b),(d),(f) of the optimal path (solution of (\ref{scaled ODE1})-(\ref{Integral condition})) for escape. (a),(b) Time profile of optimal path (green) and stable manifold $W^{s}(U_{+})$ (red) with the black marker highlighting $T_{\mathrm{end}}$ (a), $t_{\mathrm{cross}}$ (b) for a particular optimal path ($\epsilon = 1.25$, $D = 0.008$). (c),(d) Color contour plots for end time $T_{\mathrm{end}}$ (c) and the crossing time $t_{\mathrm{cross}}$ (d) in the 2 parameter $(\epsilon,D)$ - plane. (e),(f) Cross sections of (c),(d) respectively where each contour represents different value of $\epsilon$, spaced evenly at $0.04$ intervals, starting with $\epsilon = 1.05$ (dark blue, top) increasing to $\epsilon = 1.25$ (bright blue, bottom).}\label{Optimal contour}
\end{figure}

One of the advantages of reducing the study of
escape time to optimal paths is that we can perform a systematic
parameter study with moderate computational effort. First, we
investigate the timing of escape in the two parameter plane of the
ramping speed $\epsilon$ and noise level $D$, Figure \ref{Optimal contour}, panels \ref{Tend intersection}, \ref{Tend color} and \ref{Tend cross section}. Panel \ref{Tend intersection} indicates with a black marker the end time $T_{\mathrm{end}}$ at which a particular optimal path reaches the end position $x_{T} =
4$. Panel \ref{Tend color}
shows a color contour plot, with the color denoting the time
$T_\mathrm{end}$ for a range of optimal paths dependent on the ramping speed $\epsilon$ and noise level $D$. Recall that the optimal path is calculated by solving the system of
equations (\ref{scaled ODE1})-(\ref{Integral condition}) and following
the continuation steps outlined in Section~\ref{sec:steps} for a particular $\epsilon$
and large noise level $D$. Then for each $\epsilon$ a final
continuation is performed over $D$ to create an $11 \times 40$ grid
for the color plot.

The optimal path begins at $t_{0} = -10$ and so the length of the time
interval for the path is between $11.7$ and $13.5$ time units for this
range of $\epsilon$ and $D$ values. This demonstrates that for a small
$\epsilon$ value and small noise levels escape occurs for positive
$T_\mathrm{end}$, that is, with a delay compared to the time of closest
encounter of the deterministic manifolds $W^u(S_-)$ and $W^s(U_+)$
(which would be at $t=0$). As $\epsilon$ increases towards
$\epsilon_{c} = 4/3$ and the noise level increases the time to escape
decreases. This can be seen more clearly in panel \ref{Tend cross
  section}, which takes cross sections of panel \ref{Tend color} for
different values of $\epsilon$. The relationship is almost linear
between the logarithm of the noise level and the time at which the
final destination $x_{T}$ is reached. In summary, panels \ref{Tend color} and \ref{Tend cross section} indicate that the escape occurs with a delay especially for small noise. To investigate the precise value of the delay
we will look at the timing of intersection between the optimal path
and the stable manifold $W^{s}(U_{+})$, indicated by the black marker in panel \ref{tcross intersection}.

The reason for considering the stable manifold $W^{s}(U_{+})$ as a
threshold is that it plays a role similar to a saddle in stationary
escape problems. Once a realization has crossed this threshold it is
more likely to escape to $+\infty$ (in finite time in our
example). One may expect the optimal (that is, most likely) escape
path to cross this manifold when the two manifolds $W^{u}(S_{-})$ and
$W^{s}(U_{+})$ are closest together at $t = 0$. The question is then
whether the escape across the stable manifold $W^{s}(U_{+})$ occurs at
$t$ close to $0$. Panels \ref{tcross color} and \ref{tcross cross section} present the timing
of crossing the stable manifold $W^{s}(U_{+})$ to establish if this is
the case.


We observe that the range of crossing times is smaller than the range
of end times to reach $x_T$($=4$). This is expected since the
traveling time from $W^{s}(U_{+})$ to $x_{T}$ decreases for
increasing noise level $D$. For small noise levels the optimal path
tracks the manifold $W^s(U_+)$ for longer.
In the limit of large noise level $D$ in the $(\epsilon,D)$-parameter plane the most likely crossing time  is $t_{\mathrm{cross}} \approx t_{0}$
(not shown). In this limit we have a purely noise-induced transition
as the potential is nearly stationary at $t_0$. For decreasing $D$ the
intersection between the optimal path and stable manifold
$W^{s}(U_{+})$ varies with different ramping speeds $\epsilon$ such
that we have combination of noise and rate-induced tipping, with
timing depending on both parameters. As the ramping speed $\epsilon$ and noise level $D$
decreases the crossing time delay $t_\mathrm{cross}$ increases. For
the smaller noise levels in the computed range the intersection $t_{\mathrm{cross}}$ of the
optimal path with the stable manifold $W^{s}(U_{+})$ is of order
$1$, when the manifolds $W^{u}(S_{-})$ and
$W^{s}(U_{+})$ are significantly further apart than at $t = 0$. This
justifies the claim in the abstract that for noise- and rate-induced
tipping the escape is delayed in the small noise limit.
\begin{figure}[ht]
        \centering
        \subcaptionbox{Color contour.\label{M color}}[0.45\linewidth]
                {\includegraphics[scale = 0.3]{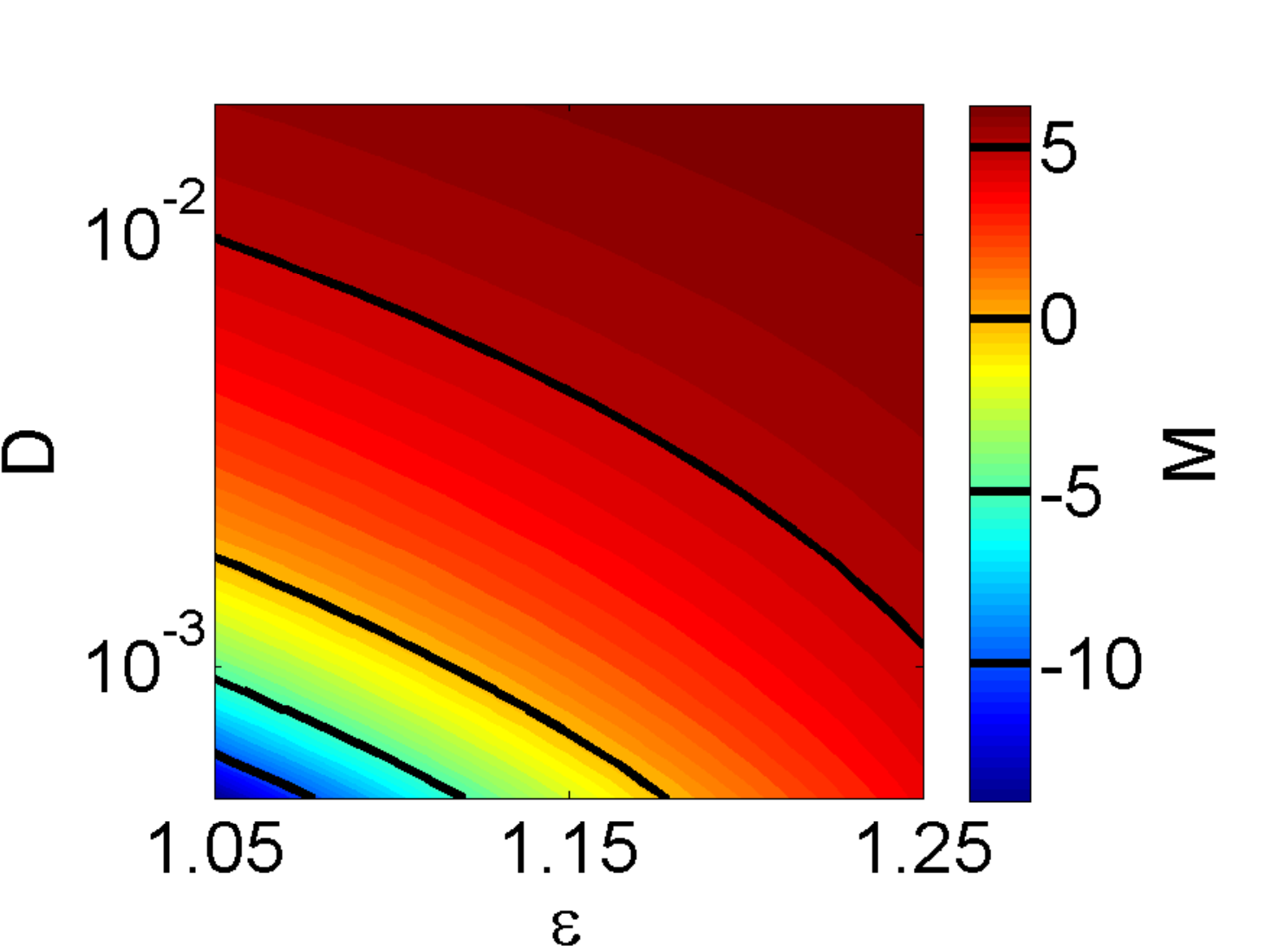}}
        \hfill 
        \subcaptionbox{Cross section.\label{M cross section}}[0.45\linewidth]
                {\includegraphics[scale = 0.3]{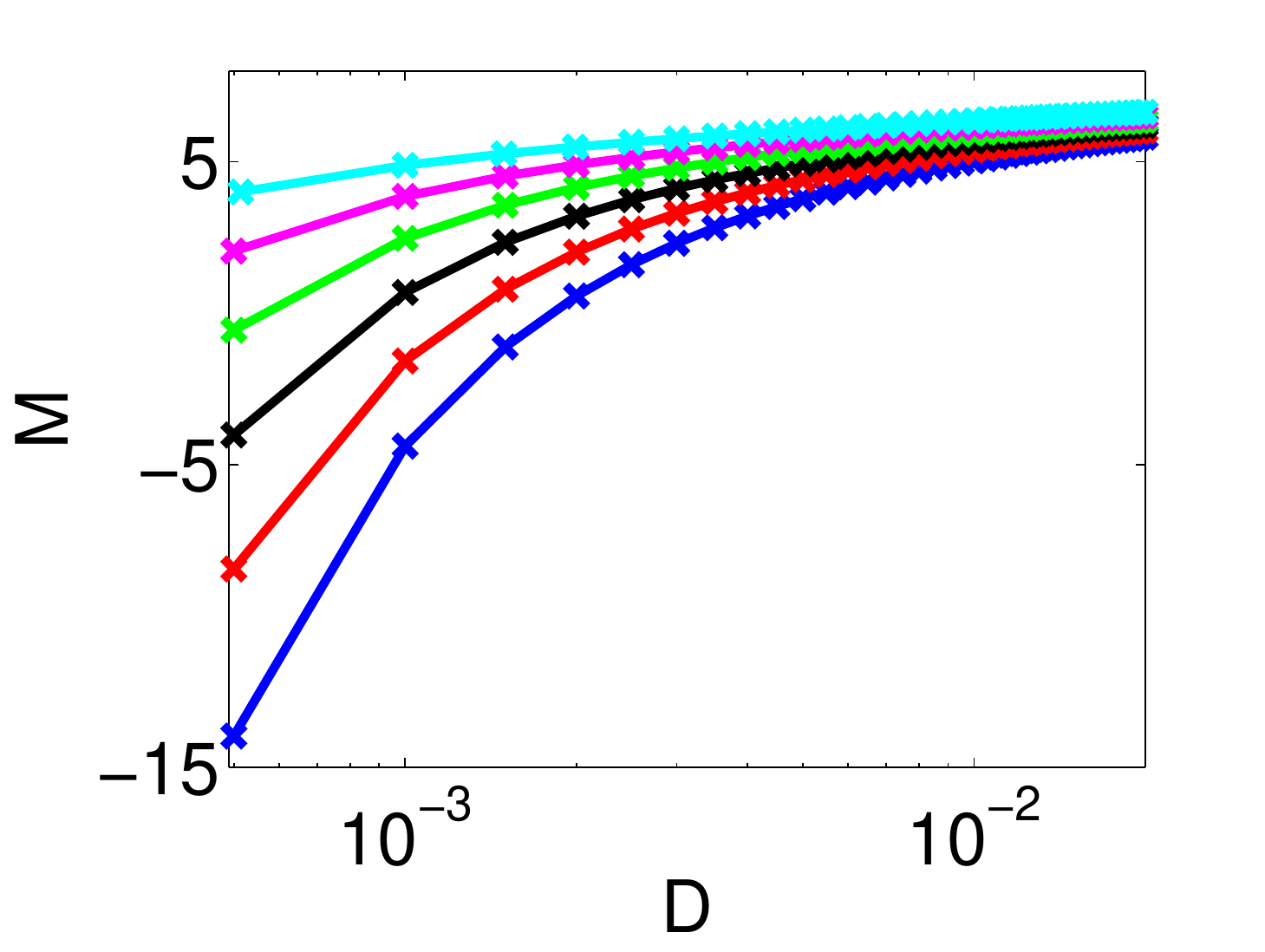}}
        ~ 
        \caption{(a) Color plot of the value of $M$. (b) Cross section of (a) each contour represents different value of $\epsilon$, spaced evenly at $0.04$ intervals, starting with $\epsilon = 1.05$ (dark blue, bottom) increasing to $\epsilon = 1.25$ (bright blue, top).}
                \label{M function plots}
\end{figure}
Figure \ref{M function plots} gives a crude estimate for the
probability of escape depending on the ramping speed $\epsilon$ and
noise level $D$. Note though, these are not the true probabilities but
rather the values of the functional $M=\log F$ (see
\eqref{Functional}) which the most likely path optimizes. The color in
Figure \ref{M color} in the 2 parameter $(\epsilon,D)$-plane equals
the value of the functional $M$ along the optimal path found at the
corresponding point in the $(\epsilon,D)$-plane. As expected, the
largest probability of escape is for large ramping speeds and large
noise levels. The value of $M$ is smallest for slow ramping speeds and
low noise levels. Figure \ref{M cross section} displaying the cross
section of Figure \ref{M color} for different values of $\epsilon$
illustrates that $M$ decreases logarithmically as $D$ decreases on a
logarithmic scale.

\section{General delay of tipping}
\label{sec:gendelay}
In Section \ref{sec:timing}, we have shown that the tipping is delayed
especially for small noise levels. 
In the context of autonomous systems \citet{bakhtin2013gumbel} gives an
asymptotic formula for this delay. \citet{bakhtin2013gumbel} considers
rare escapes for small noise levels for a process $\mathrm{d} x=xf(x)+\sqrt{2D}\mathrm{d}W_t$ on the
interval $[A,B]$ containing $0$, starting from $x_{0} < 0$, where $f$
is uniformly positive. Then the first time $T = T_{\mathrm{end}} - t_{0}$
to exit at the point $B$ (under the condition that $x(t)$ does indeed
exit at $B$) in the limit $\sqrt{2D} \rightarrow 0$  satisfies
\begin{eqnarray*}
T = c_{1}\ln\bigg(\dfrac{1}{\sqrt{2D}}\bigg)
\end{eqnarray*}
where $c_{1}$ is a constant independent of
$D$. This states that for an autonomous system, the time for rare
escapes increases linearly as the noise level decreases
exponentially. This is consistent with our findings for the
non-autonomous system, that as the noise is decreased the time
$t_{\mathrm{cross}}$ at which the optimal path crosses the stable
manifold $W^{s}(U_{+})$ increases slowly, Figure \ref{tcross cross
  section}. To conclude, we find a similar relationship for the delay
in the rate-induced tipping as that of \citet{bakhtin2013gumbel} for
rare escapes of an autonomous system. The observed level of delay in
Section~\ref{sec:timing} is of order $1$ such that the noise levels
that we consider small in Section~\ref{sec:timing} are still far
larger than the small-noise limit, for which BVPs for optimal escape
paths are available in arbitrary dimensions. (These paths tend to be
connecting orbits such that the optimal time is always infinity)\cite{ren2004minimum}.

\section{Conclusions}
\label{sec:conclusions}

We have shown that two commonly used early-warning indicators of
tipping (increase of autocorrelation and increase of variance) are
present but delayed in a prototypical model of rate-induced
tipping. By looking at the timing of escape using optimal paths we
find that the tipping event itself is delayed for small noise
levels. We conclude that the delay in the early-warning indicators is
consistent with the delay in the actual tipping (at least for the
example).

We extended the boundary-value problem for the most likely path for
tipping (escape) based on \citet{zhang2008theory} to include
optimality of time for finite noise. This additional optimality
criterion created a variational optimization problem that we solved
computationally with continuation techniques (using the package
AUTO). With the help of continuation we performed a systematic
parameter study in the $(\epsilon,D)$-plane (ramping speed vs. noise
level). The time when the optimal path for escape crosses the stable
manifold $W^{s}(U_{+})$ is a measure for the timing of tipping. We
find that for large ramping speeds and noise levels there is no delay
and even for lower ramping speeds there is only a small
delay. However, for small noise levels $D$ the tipping delay is of
order $1$. 

\change{We hypothesize that the observed delay in tipping is present
  \change{independent of the particular form of $\lambda(t)$ as long as it
    is qualitatively similar to the ramp like shift \eqref{lambda prime}. Similarly, this delay should be observable independent of
    the particular shape of the potential well $U(\cdot,t)$.}} This
paper demonstrated that the optimal path for escape, a solution of a
BVP, matches simulation results well. \change{The technique used to
  find the optimal path of escape finds the local maximum and is
  general such that it can be used to determine the timing for any
  type of tipping.}

However, \change{the optimal path} may miss
the global optimum when there is more than one realistic opportunity
for escape. For a small single window of escape as considered in this
paper the escape rate will form a unimodal distribution with a narrow
peak, for which the optimal path is close to the mode. However, if one
considers different scenarios $\lambda(t)$ for ramping the system
parameter (for example, one that is not monotonically increasing, see
\citet{ashwin2015parameter}), the escape rate would have a multimodal
distribution. We conjecture that we find one optimal path for each of
the modes of the distribution. It is unclear if for non-monotone
parameter shift $\lambda(t)$ the tipping or the early-warning
indicators are delayed for the small noise levels.  This would further
support the conclusion that the autocorrelation and variance can be
used as early-warning signals for rate-induced tipping events. Furthermore,
this paper has focused on the one-dimensional case. Thus, an extension
to the general multiple dimensional case is still required.



%
%

%



\begin{appendix}
\setcounter{paragraph}{0}
\section{\change{Dependence on parameters}}\label{sec:threshold}
\change{ This appendix details how the choice of the parameters
  $x_{\mathrm{end}}$, the upper domain boundary, and the threshold
  parameter $y$ affects the results presented in the paper.}
\subsection{\change{Domain boundary parameter}}
\change{We investigate the effect the upper boundary of the domain,
  $x_{\mathrm{end}}$, has on the early-warning indicators, the
  increase of autocorrelation and variance. We choose a domain
  \change{$[x_\mathrm{start},x_\mathrm{end}]$ with $x_\mathrm{start}=-6$ and
    $x_\mathrm{end}$ as shown in Figure~\ref{EWI analysis}.} Choosing
  the domain fixed \change{in time} is natural as in realistic problems
  \change{we do not know a-priori the location of the moving well}.  In
  Section~\ref{sec:apparent failure of early-warning} we chose
  $x_{\mathrm{end}} = 2$, which corresponds to a wide domain for the
  problem \change{(including realizations into the computation of the
    early-warning indicator, which are already escaping)}. Figure
  \ref{EWI analysis} explores the effect narrowing the domain has on
  the decay rate estimate \change{$\theta=(1-a)/\Delta t$ (where $a$ is the lag-$1$
  autocorrelation with time step $\Delta t$)} and variance $V$. We
  consider the decay rate, instead of the linearly related
  autocorrelation as the decay rate is independent of the time step
  $\Delta t$.

\begin{figure}[ht]
        \centering
        \subcaptionbox{Decay rate\label{Autocorrelation analysis}}[0.45\linewidth]
                {\includegraphics[scale = 0.3]{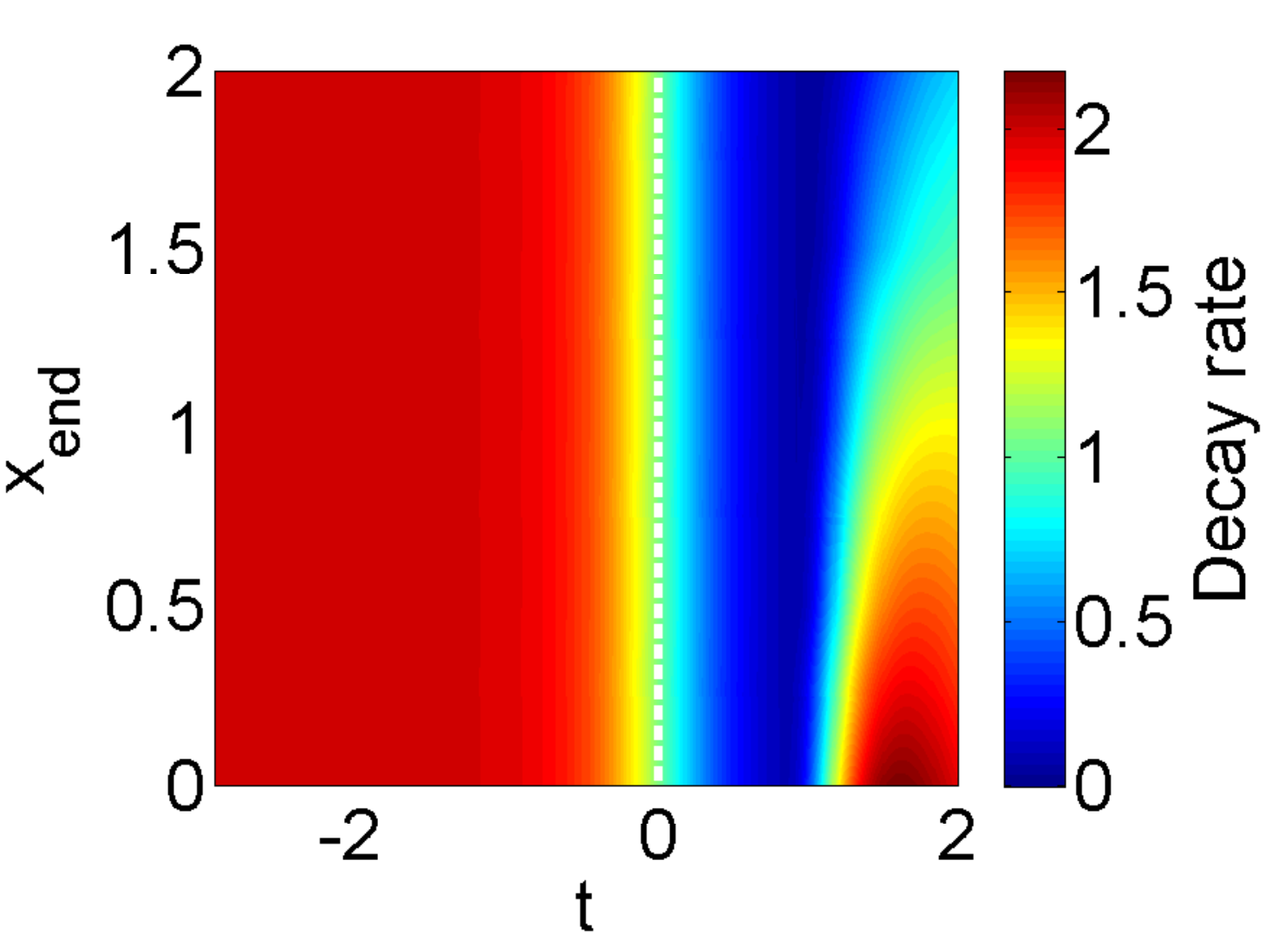}}
        \hfill 
        \subcaptionbox{Variance\label{Variance analysis}}[0.45\linewidth]
                {\includegraphics[scale = 0.3]{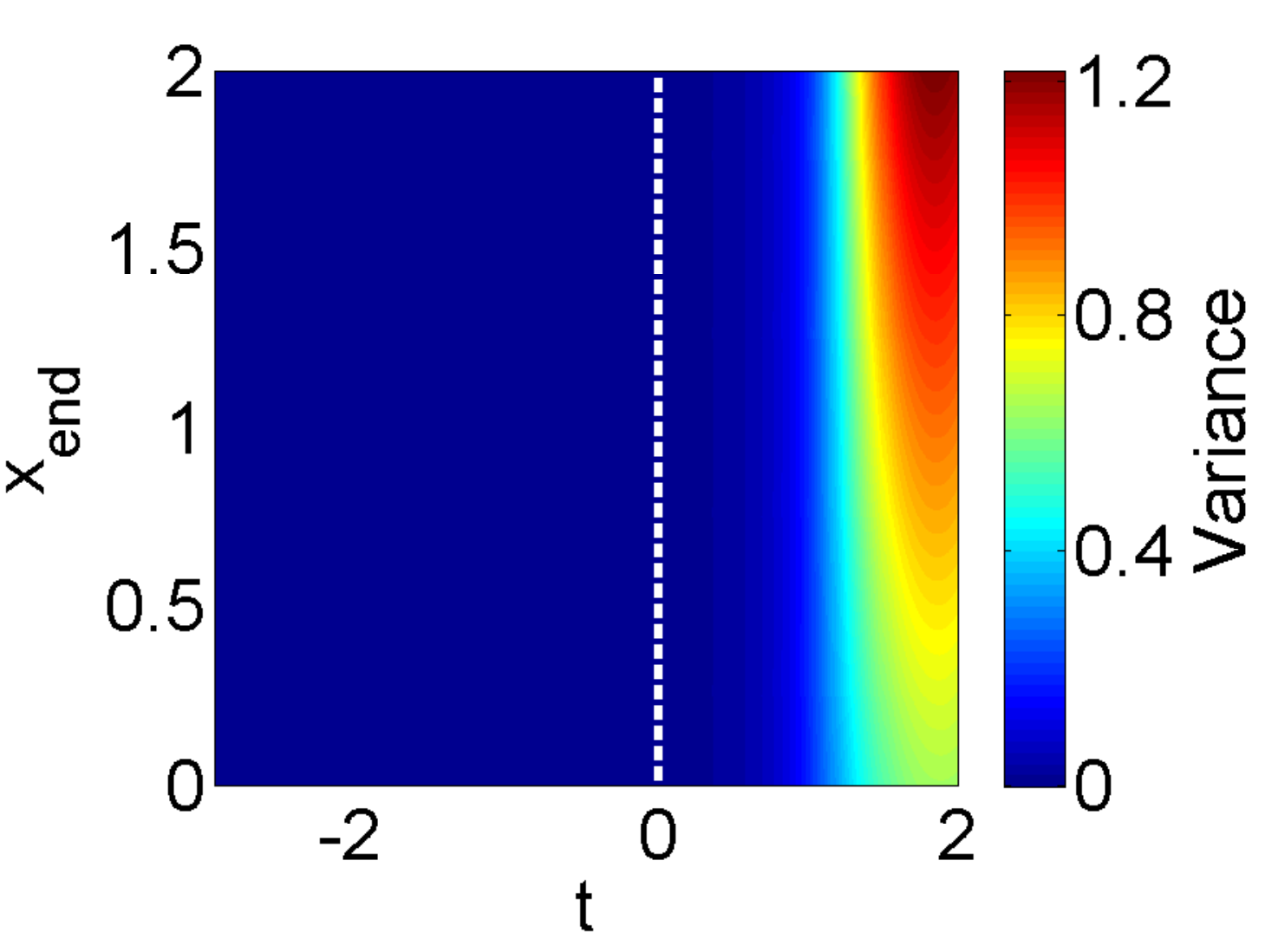}}
        ~ 
        \caption{\change{Effect of the width of domain on the decay rate and variance by varying upper boundary $x_{\mathrm{end}}$}}\label{EWI analysis}
\end{figure}

In Figure \ref{Autocorrelation analysis}, we see that the \change{timing
  of the} onset of the \change{decrease of the} decay rate estimate
$\theta$ is independent of the upper boundary of the domain,
$x_{\mathrm{end}}$, and hence, so is \change{timing of the onset of} the \change{increase of the}
lag-$1$ autocorrelation $a$. Likewise, in Figure \ref{Variance
  analysis} the \change{timing of the} onset of the variance is
independent of $x_{\mathrm{end}}$ and, importantly, shows no increase
before $t = 0$. The precise values of the autocorrelation and the
variance depend, of course, strongly on the width of the domain (and,
hence, on $x_\mathrm{end}$). Thus, we have shown that, while
autocorrelation and variance change quantitatively with the domain
width, the timing of their increase (which is the early-warning
indicator) does not change.}

\subsection{\change{Threshold parameter}}\label{subapp:threshold}
\begin{figure}[ht]
        \centering
        \subcaptionbox{$D = 0.1$\label{Escape rate color large D}}[0.45\linewidth]
                {\includegraphics[scale = 0.3]{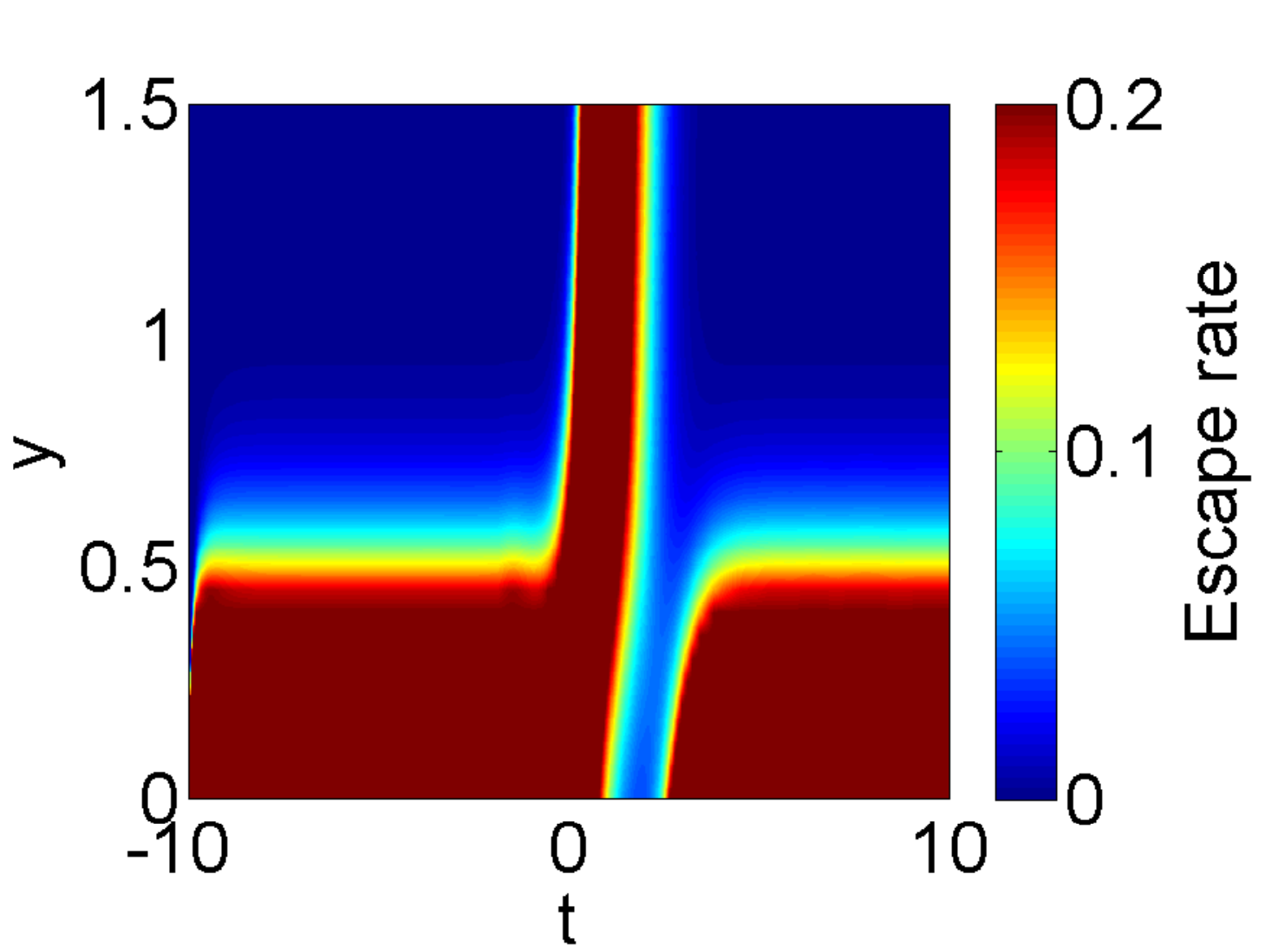}}
        \hfill 
        \subcaptionbox{$D = 0.008$\label{Escape rate color small D}}[0.45\linewidth]
                {\includegraphics[scale = 0.3]{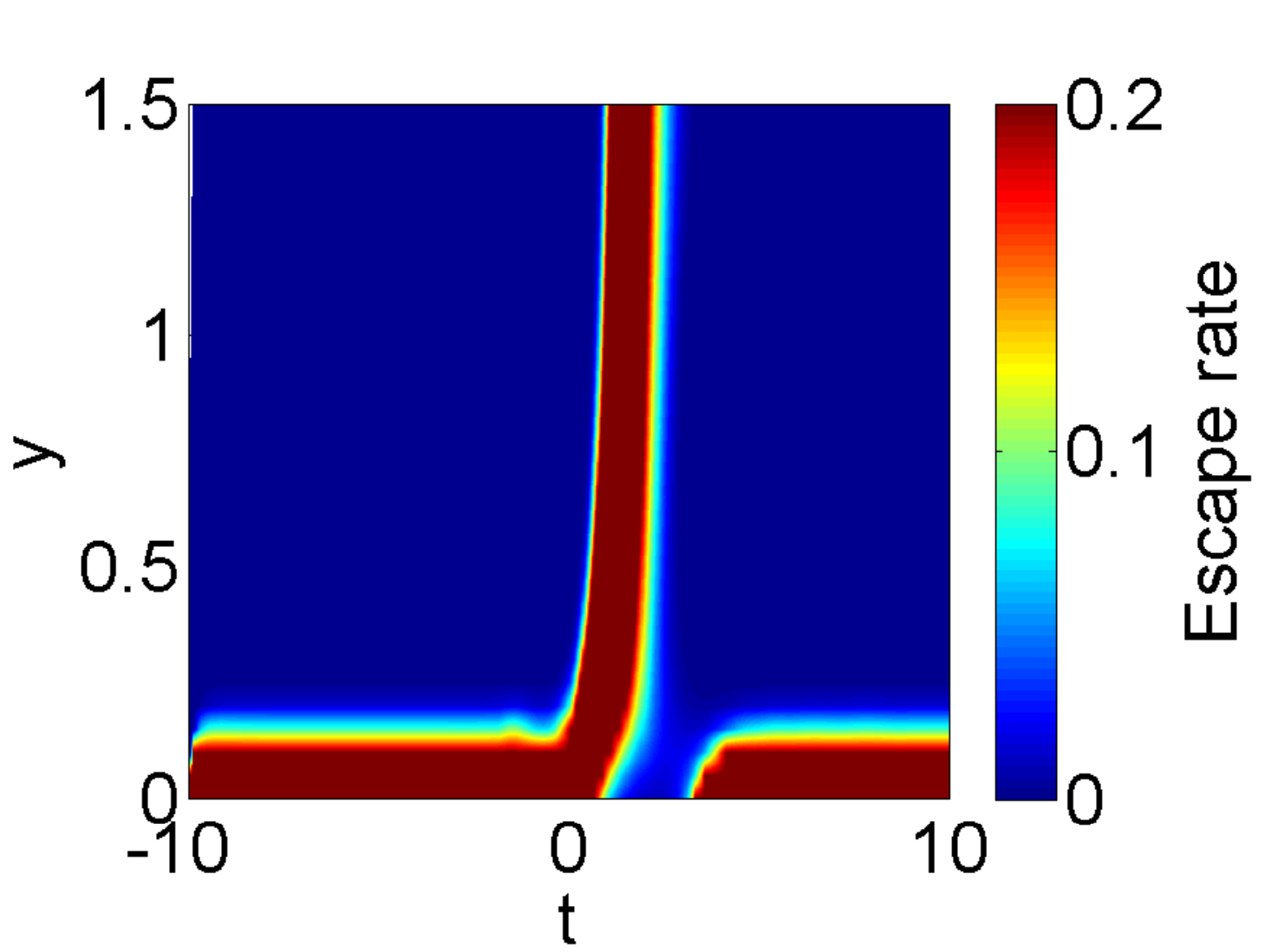}}
        ~ 
        \caption{Evaluating distance $y$ required between the deterministic trajectory $x^{u}(t)$ and threshold $\tilde{x}(t)$ in both the large and small noise limit cases, $\epsilon = 1.25$}\label{Escape rate color}
\end{figure}

This \change{section} presents in more detail how the distance $y$ of the
threshold curve $\tilde{x}(t)$ (at which we consider a realization as
escaped) from the deterministic trajectory $x^u(t)$ influences our
results ($\tilde{x}(t)= x^{u}(t) + y$). If we choose $y$ too small, then escape will be detected everywhere. This is demonstrated by Figure \ref{Escape rate color}. For example, when the steady state is stationary no escape should be detected, because, if realizations cross the threshold $\tilde{x}(t)$, most will not escape to $+\infty$ but will return to the unstable manifold $W^{u}(S_{-})$. Clearly for a larger noise level $D$ a greater $y$ is required as there will be larger fluctuations about the unstable manifold $W^{u}(S_{-})$ than for a small noise level.

Figure \ref{Escape rate color large D} demonstrates that for values
of $y$ less than a critical value $y_{c} \approx 1$, a large fraction
of realizations would count as escaped even for $t$ close to $\pm 10$,
where $\lambda$ is close to stationary, in the large noise case. For
smaller noise level Figure \ref{Escape rate color small D}
illustrates that a smaller $y$ will suffice. Both figures show that
above a certain minimal value of $y$ the window of escape remains
nearly independent of $y$. The value $y=1.5$, used in
Section~\ref{sec:optimal} for our threshold curve $\tilde{x}(t)$, is
well above that minimal value of $y$.

\section{Variational optimization problem for specific example}\label{sec:bvpderiv}
The following is the variational optimization problem for the specific rate-induced example discussed in the paper, equation \eqref{SDE}.

The potential $U(x,\lambda(t))$ for equation \eqref{SDE} is:
\begin{eqnarray*}
U(x,\lambda(t)) &&= -\dfrac{x^{3}}{3} - \lambda x^{2} + (1-\lambda^{2})x \mbox{\quad such that}\\
U' &&= -x^{2} - 2\lambda x + 1 - \lambda^{2} \\
U'' &&= -2(x+\lambda) \\
\dot{U} &&= -\dot{\lambda}x(x + 2\lambda) = -\epsilon\lambda x(\lambda_{\max} - \lambda)(x + 2\lambda)
\end{eqnarray*}
where $U'$ and $\dot{U}$ represent the derivatives of the potential w.r.t. space and time respectively and $\dot{\lambda}$ is given by equation \eqref{rtip lambda dot init}. These equations feed into the $V_{s}$, equation \eqref{Vs equation} giving:

\begin{eqnarray*}
V_{s} = &&\dfrac{x^{4} + 4\lambda^{2}x^{2} + (1-\lambda^{2})^{2} + 4\lambda x^{3} - 2x(x + 2\lambda)(1-\lambda^{2})}{4D} \\
&&+ x + \lambda + \dfrac{\epsilon\lambda x(\lambda_{\max} - \lambda)(x + 2\lambda)}{2D}
\end{eqnarray*}
and so the 2nd order boundary value problem, equation \eqref{2nd order ODE}, split into 2 first order ODEs augmented with the ODE for $\lambda$, equation (\ref{rtip lambda dot init}), which are to be solved on the $[0,1]$ time domain in AUTO looks like:

\begin{eqnarray}
\label{rtip x1 dot}
\dot{x}_{1} &&= x_{2}(T_{\mathrm{end}} - t_{0})\\
\label{rtip x2 dot}
\dot{x}_{2} &&= h_{2}(x_{1},\lambda(t))(T_{\mathrm{end}} - t_{0})\\
\label{rtip lambda dot}
\dot{\lambda} &&= h_{3}(\lambda(t))(T_{\mathrm{end}} - t_{0})
\end{eqnarray} 
where \eqref{rtip x1 dot}-\eqref{rtip lambda dot} correspond with \eqref{scaled ODE1}-\eqref{scaled ODE2} and

\begin{eqnarray*}
h_{2}(x_{1},\lambda(t)) &&= 2D\dfrac{\partial V_{s}(x_{1},\lambda(t))}{\partial x_{1}} \\
h_{3}(\lambda(t)) &&= \epsilon\lambda(\lambda_{\max} - \lambda)
\end{eqnarray*}
The function $M = \log(F)$, equation \eqref{Integral M} for the general case, is maximized and used to monitor any maxima or minima, is given by:

\begin{eqnarray}
\nonumber
M = &&\int_{0}^{1}\bigg[\dfrac{U(x_{0},\lambda(t_{0})) - U(x_{T},\lambda(T_{\mathrm{end}}))}{2D} \\
&&- \bigg(\dfrac{x_{2}}{4D} + V_{s}(x_{1},\lambda(t))\bigg)(T_{\mathrm{end}}-t_{0})\bigg]\mathrm{d}t
\label{Integral condition 1}
\end{eqnarray}
The variational equations for $z_{1}$, $z_{2}$ \eqref{z ODE} and $z_{3} = \dfrac{\partial \lambda(t)}{\partial T_{\mathrm{end}}}$ are given as:

\begin{eqnarray}
\label{rtip z1 dot}
\dot{z}_{1} = &&x_{2} + z_{2}(T_{\mathrm{end}} - t_{0}) \\
\label{rtip z2 dot}
\dot{z}_{2} = &&h_{2}(x_{1},\lambda(t)) \\
\nonumber
&&+ \bigg(\dfrac{\partial h_{2}(x_{1},\lambda(t))}{\partial x_{1}}z_{1} + \dfrac{\partial h_{2}(x_{1},\lambda(t))}{\partial \lambda}z_{3}\bigg)(T_{\mathrm{end}} - t_{0}) \\
\label{rtip z3 dot}
\dot{z}_{3} = &&h_{3}(\lambda(t)) + \dfrac{\mathrm{d}h_{3}(\lambda(t))}{\mathrm{d}\lambda}(T_{\mathrm{end}} - t_{0})
\end{eqnarray}

To locate the local maximum of $M$, which is the derivative of \eqref{Integral condition 1} w.r.t. $T_{\mathrm{end}}$ we have a second integral condition corresponding to equation \eqref{Integral condition} in the paper (multiplied by $-4D$ to remove fractions): 

\begin{eqnarray}
\nonumber
&&\int_{0}^{1}\bigg[\bigg(z_{2} + 2h_{2}(x_{1},\lambda(t))z_{1} + \dfrac{\partial V_{s}(x_{1},\lambda(t))}{\partial\lambda}z_{3}\bigg)(T_{\mathrm{end}} - t_{0}) \\ &&+2\dfrac{\partial U(x_{T},\lambda(T_{\mathrm{end}}))}{\partial T_{\mathrm{end}}} + x_{2} + 4DV_{s}(x_1,\lambda(t))\bigg]\mathrm{d}t = 0
\label{Integral condition 2}
\end{eqnarray}
and thus for the general example in the paper which had five equations to solve \eqref{scaled ODE1}-\eqref{scaled ODE2}, \eqref{z ODE} and \eqref{Integral condition} and for the specific rate-induced example there are seven equations to solve \eqref{rtip x1 dot}-\eqref{rtip lambda dot}, \eqref{rtip z1 dot}-\eqref{rtip z3 dot} and \eqref{Integral condition 2}.

\section{Detailed explanation of continuation steps}\label{app:continuation steps}

We provide further explanation of the continuation steps presented in the paper for the specific rate-induced example \eqref{rtip x1 dot}-\eqref{rtip lambda dot}.

\subsection{Step 1: $T_{\mathrm{init}}$ continuation} This
continuation is similar to an integration in time
continuation. However, there is a difference between this continuation
and performing a continuation in $T_{\mathrm{end}}$ to $-9$ having
started with $T_{\mathrm{end}} = -10$ and $T_{\mathrm{init}} =
1$. Starting with $T_{\mathrm{end}} = -10$ would mean:

\begin{eqnarray*}
\dot{x}_{1} &= 0 \qquad x_{1}(-10) &= -1 \\
\dot{x}_{2} &= 0 \qquad x_{1}(-9) &= -1
\end{eqnarray*}
which has no unique solution. Incorporating the artificial continuation parameter $T_{\mathrm{init}}$ (initially at $0$) and setting $T_{\mathrm{end}} = -9$ we have:

\begin{eqnarray*}
\dot{x}_{1} &= x_{2} \quad x_{1}(-10) &= -1 \\
\dot{x}_{2} &= 0 \quad x_{1}(-9) &= -1
\end{eqnarray*}
which does have a locally unique solution. Thus, we can continue in
$T_\mathrm{init}$ until $T_{\mathrm{init}} = 1$ to obtain a solution
of the full system \eqref{rtip x1 dot}--\eqref{rtip lambda
  dot}. Note that the parameters $M$ and $m$ have to be kept free
during this continuation such that the integral conditions
\eqref{Integral condition 1}, \eqref{Integral condition 2} are always
satisfied.

\subsection{Step 2: $x_{T}$ continuation}\label{app:xT continuation} The sketch of the phase
portrait, Figure \ref{PP sketch}, is an accurate representation of the
full system that does not change over the time considered, since this
continuation is for $t_{0} = -10$ and $T_{\mathrm{end}} = -9$ and
hence $\lambda \approx 0$. The phase portrait contains two saddles,
which are located close to the equilibrium points of $S_{-}$ and
$U_{-}$ and one center close to the origin. The saddles are offset to
the left of $S_{-}$ and $U_{-}$ by approximately $D/2$. The center of
the elliptic region is located at $(x_{1},x_{2}) \approx (D,0)$.

The initial trajectory after step~1 ($x_{T} = x_{0} = -1$) is a unique
solution that is contained within the elliptic region near the saddle
(see blue trajectory Figure \ref{PP sketch}). The boundary condition
$x_{1}(t_{0}) = x_{0}$ ensures that all trajectories start on the
dashed line $S_{-} = -1$. During continuation in $x_T$ the
trajectories need to travel further but still in the same time
interval ($t_{0} = -10$, $T_{\mathrm{end}} = -9$) as $x_{T}$ is
increased. Therefore, the starting position increases in the $x_{2}$
direction where the vector field has a larger $x_1$ component. This
enables the trajectories to travel faster in order to travel further
in the same time period. See Figure \ref{PP sketch} for intermediate
phase portraits in this continuation, Figure \ref{PP after xT
  continuation} for the phase portrait and \ref{Optimal path after xT
  continuation} for the time profile of the (final) trajectory with
$x_T=4$.

\begin{figure}[ht]
        \centering
        \subcaptionbox{$x_{T} = 4$\label{PP after xT continuation}}[0.45\linewidth]
                {\includegraphics[scale = 0.3]{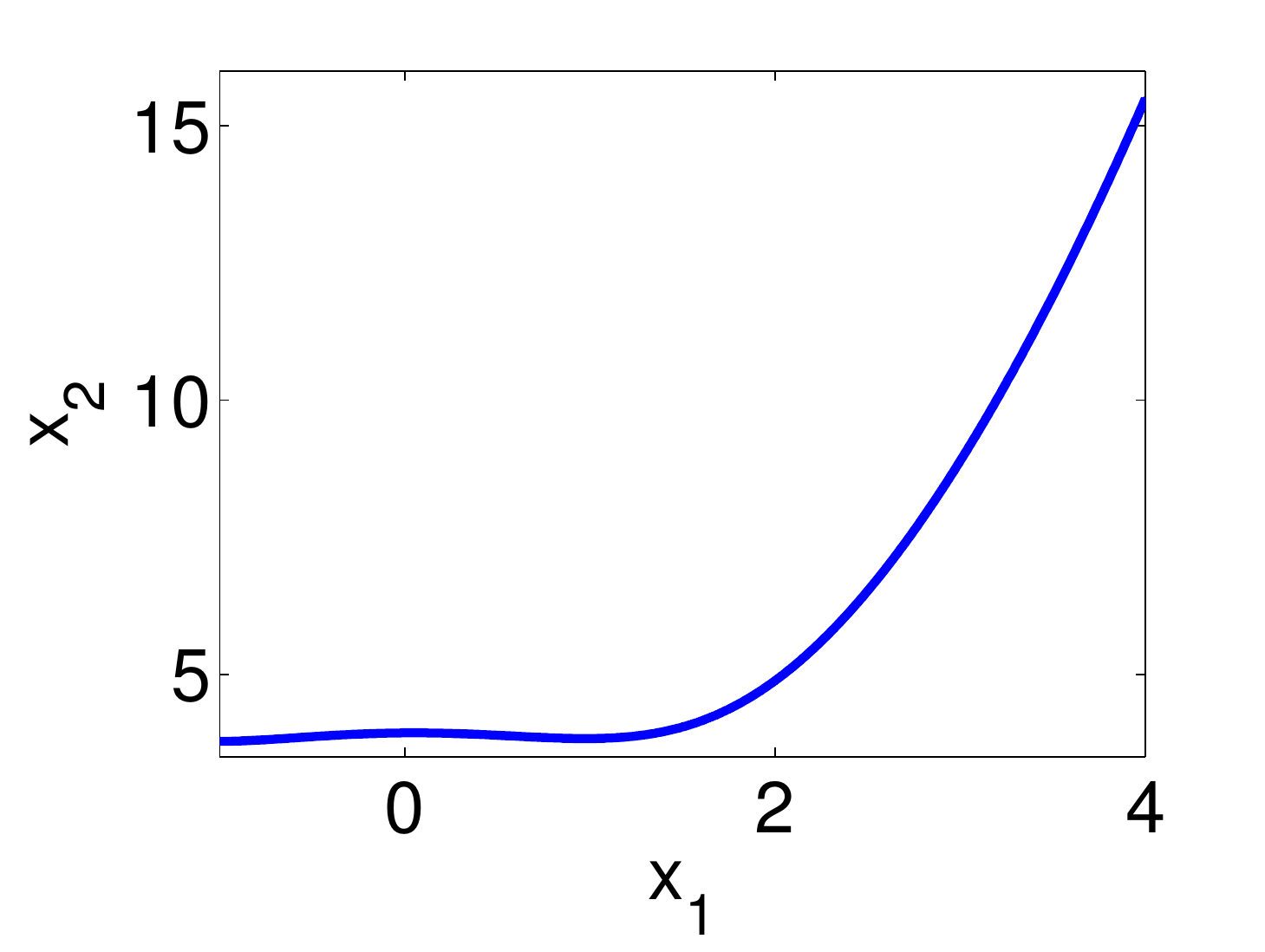}}
        \hfill 
        \subcaptionbox{Optimal path after $x_{T}$ continuation step.\label{Optimal path after xT continuation}}[0.45\linewidth]
                {\includegraphics[scale = 0.3]{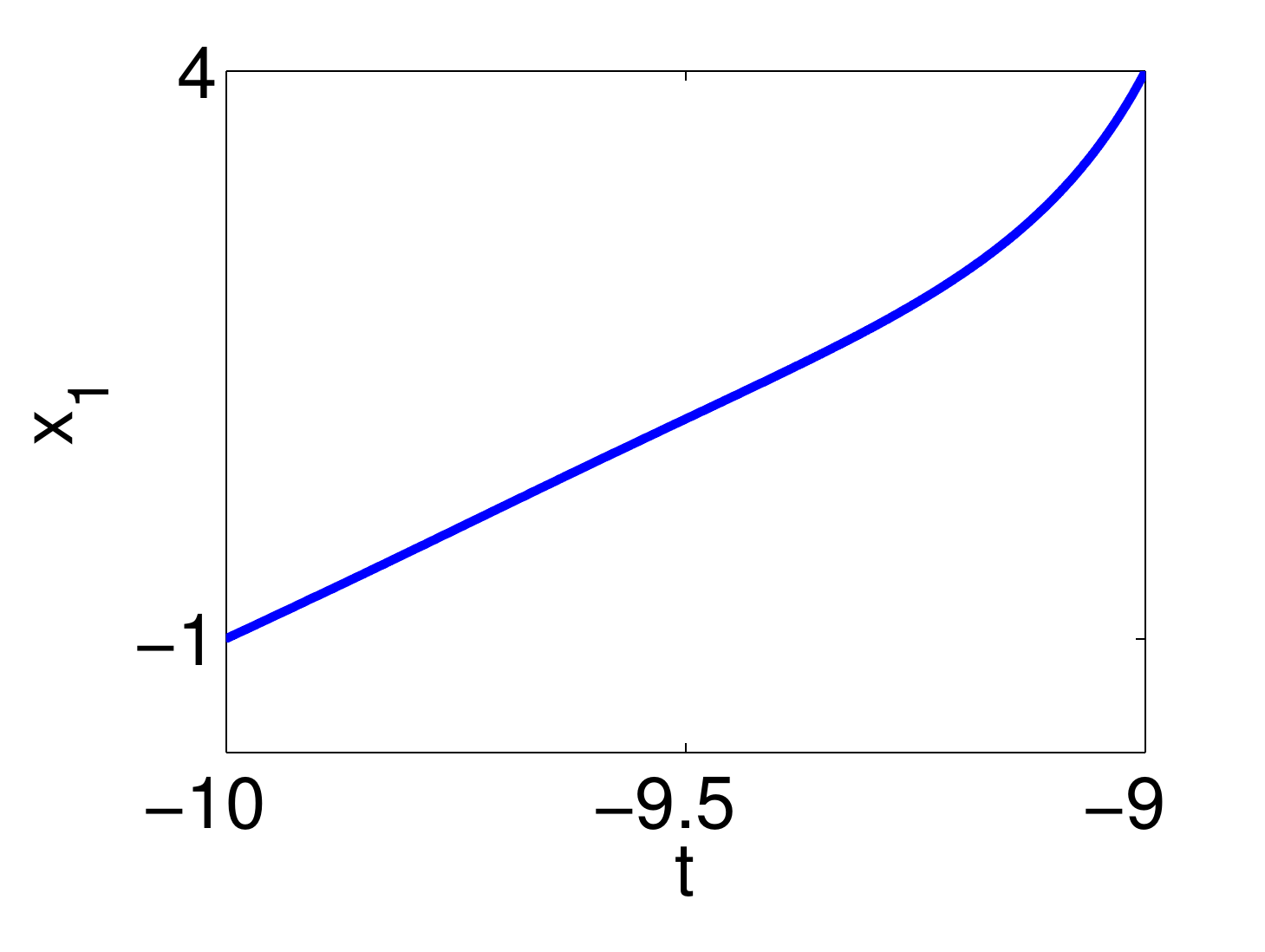}}
        ~ 
        \caption{(a) Trajectory in $(x_{1},x_{2})$ - plane, (b) optimal path after $x_{T}$ continuation. $t_{0} = -10$, $\epsilon = 1.25$, $D = 0.05$.}\label{End of xt continuation}
\end{figure}
This trajectory corresponds to the optimal path for a purely
noise-induced tipping since $\lambda$ is close to stationary. One
would expect the optimal time for tipping to be a result of both noise
and rate-induced tipping.
   
\subsection{Step 3: $T_{\mathrm{end}}$ continuation} The next step is
to perform a continuation in $T_{\mathrm{end}}$ while monitoring $m$
for roots (or $M$ for critical points). Since $M$ may have several
local minima and maxima for increasing $T_\mathrm{end}$
we continue in $T_{\mathrm{end}}$ too sufficiently large values where
we observe the asymptotic monotone decrease of $M$. 
In our example, we continued $T_{\mathrm{end}}$ from
$-9$ to $20$, monitoring the bifurcation diagram in
the $(T_{\mathrm{end}},M)$ - plane (not shown but similar to Figure
\ref{Tend continuation 1}).

\begin{figure}[ht]
        \centering
        \subcaptionbox{$T_{\mathrm{end}} = 20$.\label{opt path Tend20}}[0.45\linewidth]
                {\includegraphics[scale = 0.3]{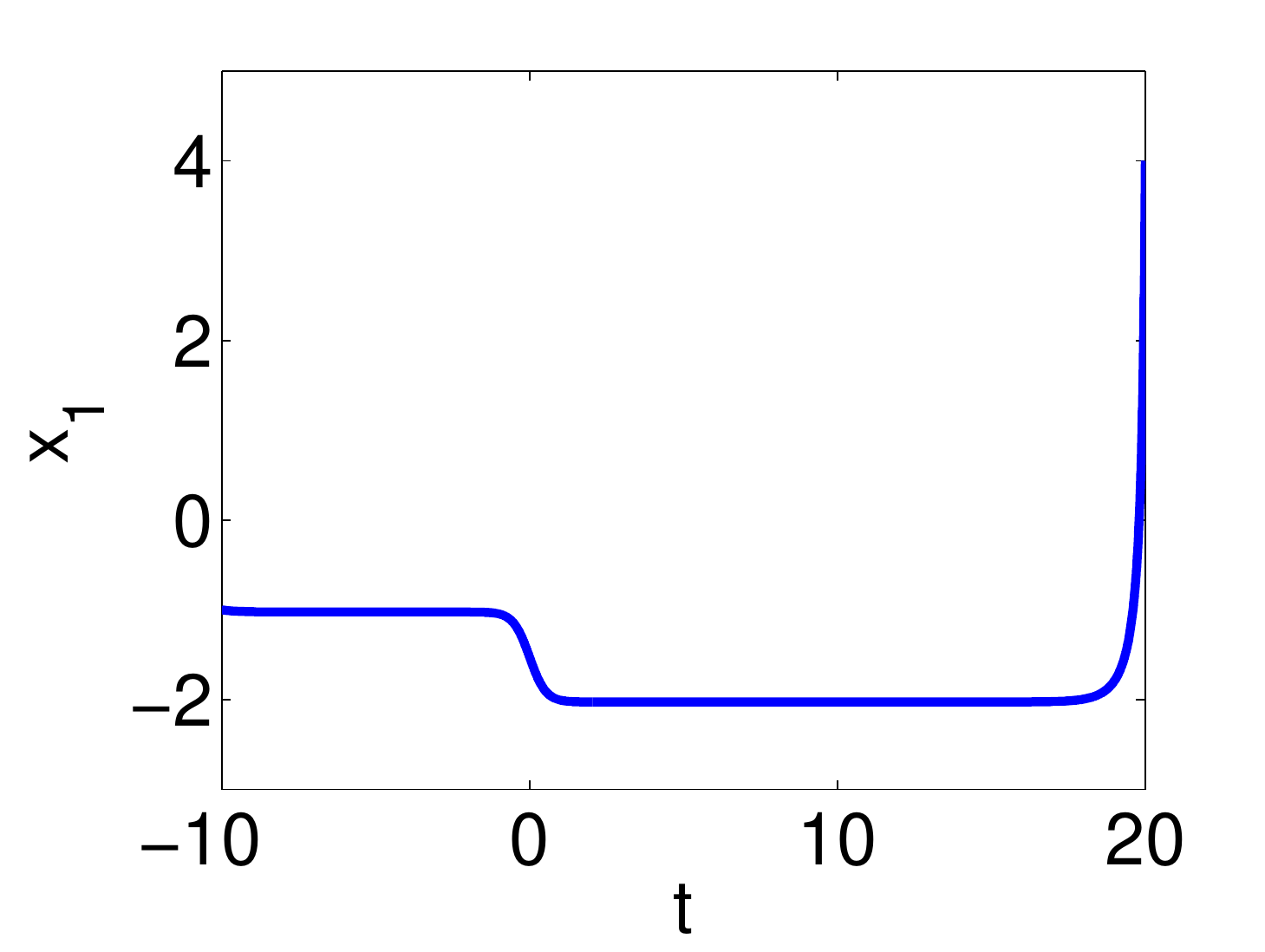}}
        \hfill 
        \subcaptionbox{$T_{\mathrm{end}} \approx 1.43$.\label{dM continuation}}[0.45\linewidth]
                {\includegraphics[scale = 0.3]{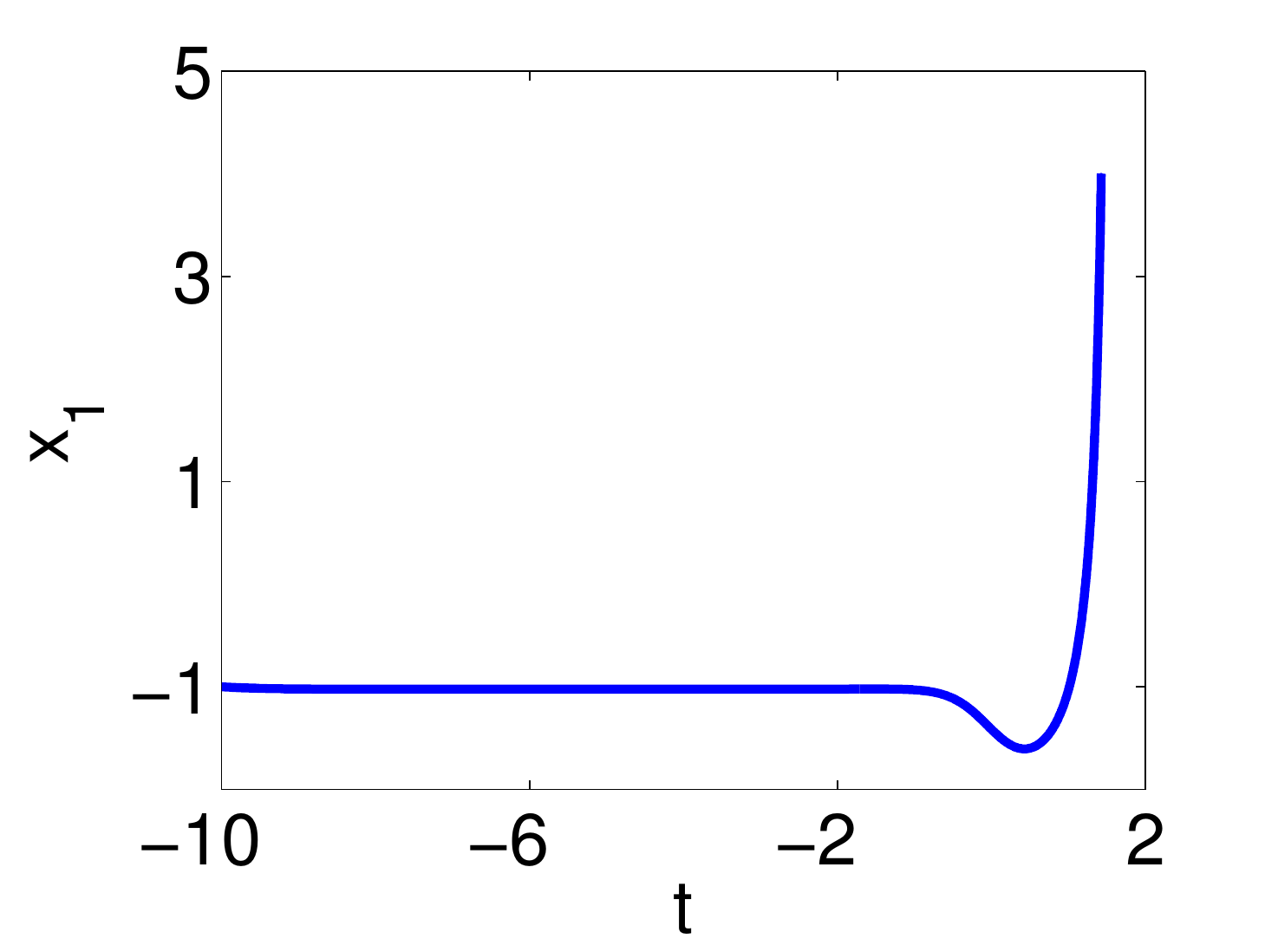}}
        ~ 
        \caption{Comparison between optimal paths for $T_{\mathrm{end}} = 20$ and after $m$ continuation is completed. $t_{0} = -10$, $\epsilon = 1.25$, $D = 0.05$}\label{dM and D continuations}
\end{figure}

For choice of parameter values displayed in Figure~\ref{Tend
  continuation 1} there is only one critical point of $M$, which
corresponds to the maximum we are interested in. Figure \ref{opt path
  Tend20} displays the optimal path after the $T_{\mathrm{end}}$
continuation step, for $T_{\mathrm{end}} = 20$. For $T_{\mathrm{end}} = 20$ the optimal path reaches the saddle at $t \approx 1$, but waits at the saddle until $t \approx 18$ before escaping to $x_{T}$, which \change{is optimal only in the limit $D\to0$}. Therefore, we detect when $m = 0$ which satisfies the integral condition \eqref{Integral condition 2} to identify the maximum of $M$ and hence achieve the optimal path in an optimal time, see Figure \ref{dM continuation}.


\end{appendix}

\bibliography{Rate-induced_paper}

\end{document}